\documentclass{article}
\usepackage{amsmath,amssymb,bm,yhmath,mathtools,extarrows,MnSymbol}
\usepackage[all]{xy}
\usepackage{caption,graphicx}
\font\ibf=cmbxti10

\usepackage{multicol}

\title{Topics in geometric group theory, Part II}
\author{by Valentin {\sc Po\'enaru}\footnote{Professor Emeritus at the Universit\'e Paris Sud-Orsay, Math\'ematiques 425, 91405 Orsay Cedex, France.  \newline e-mail: valpoe@hotmail.com}}
\date{\sl (November 2016, revised November 2017)}

\begin{document}

\maketitle

\begin{abstract}
This paper gives a quick overview of the author's recent result that all finitely presented groups are QSF.
\end{abstract}

\setcounter{section}{-1}
\section{Introduction}\label{sec1}

In this second part of our survey we discuss the recent result that all finitely presented groups have the QSF property, recently proved by the author. We give here an outline of the complete long proof, contained in the trilogy of papers \cite{21}, \cite{22}, \cite{23}. We hope that the present outline should help in understanding the much longer trilogy.

\smallskip

The QSF concept, due to Brick, Mihalik and Stallings is reviewed in the beginning of the next section. It has its roots in the old seminal work of Max Dehn and as the reader will be able to see, it is not by chance that our two parts survey is dedicated to the memory of Dehn (and more about Dehn can be found in the paper of Pierre de la Harpe ``Topologie, th\'eorie des groupes et probl\`emes de d\'ecision'', Gazette des Math\'ematiciens, n$^{\rm o}$~125 (2010) pp. 41-75).

\smallskip

Part I of this survey should, among other things, give the necessary preliminary material for what is being told here, and also give the context and motivation for the QSF theorem.

\smallskip

The present survey paper is divided into three parts, the third one being a road-map for the whole approach. The reader should not feel obliged to wait untill the end, to look at it.

\smallskip

I wish to thank C\'ecile Gourgues for the typing and Marie-Claude Vergne for the drawings. And I also want to thank the anonymous referee who helped me write a better paper.

\section{QSF Theory}

I will start by explaining the QSF notion for finitely presented groups. So I will review now some important work of S. Brick, M. Mihalik and the late J. Stallings \cite{1}, \cite{2}.

\smallskip

They introduced the following notion.

\bigskip

\noindent {\bf Definition (1.1) --} A locally compact complex $X$ (and from now on we will stay always inside the simplicial category) is said to be QSF if for every compact $k \subset X$ there is a finite simply-connected complex $K$ endowed with a simplicial map $f$ into $X$, and also with an injection $k \xrightarrow{ \ j \ } K$, all entering in the commutative diagram below
$$
\xymatrix{
k \ar[rr]_i \ar[rd]^j &&X \, , \\
&K \ar[ur]_f
}
\leqno (1.2)
$$
where the following Dehn-like condition is also fulfilled
$$
jk \cap M_2 (f) = \emptyset \, .
$$

This looks, of course, very much like the definition of Dehn-exhaustibility from the first part of this paper, after which it might well have been modelled, but in (1.2) the $f$ is just a simplicial map and not necessarily an immersion.

\smallskip

The ``QSF'' comes from ``quasi-simply-filtered'', but I am not responsible for this terminology. End of definition (1.1).

\bigskip

Here comes now one of the big main virtues of (1.1).

\bigskip

\noindent {\bf Lemma 1.} (Brick, Mihalik and Stallings) {\it Let $G$ be a finitely presented group and let also $P_1 , P_2$ be two compact complexes such that $\pi_1 P_2 = \pi_1 P_2 = G$, i.e. two geometric presentations of $G$. Then $\widetilde P_2$ is QSF iff $\widetilde P_1$ is QSF and in this case, one says that the group $G$ is QSF.}

\bigskip

So QSF is a group-theoretical notion, unlike Dehn-exhaustibility, or WGSC, or GSC. There is also the following fact, analogous to the Dehn-lemma \`a la Po, from Part I of this paper, namely

\bigskip

\noindent (1.3) \quad If $M^3$ is a closed 3-manifold and if $\pi_1 M^3$ in QSF, then $\pi_1^{\infty} \widetilde M^3 = 0$.

\bigskip

For these basic elementary facts, the references are \cite{1}, \cite{2}. Retain, anyway that it makes sense to talk about a group $G$ as being (or not being) QSF and of $\pi_1^{\infty} G$ too. But remember also that via the work of M.~Davis \cite{15} not every finitely presented group $G$ comes with $\pi_1^{\infty} G = 0$, even if $K(G,1)$ is a compact manifold.

\smallskip

Now, it follows from the celebrated work of Grisha Perelman on the 3-dimensional Poincar\'e Conjecture and the Geometrization of 3-manifolds \cite{5}, \cite{6}, \cite{7}, \cite{4}, \cite{8}, \cite{39}, and see also \cite{16}, \cite{17}, that we have the following result (which is also implied by the more recent work of Agol-Wise \cite{53}, \cite{54}; and see here, also, the very comprehensive reference \cite{52}).

\bigskip

\noindent {\bf Theorem 2.} (G. Perelman) {\it If $M^3$ is a closed $3$-manifold, then $\pi_1 M^3$ is QSF, and hence, also $\pi_1^{\infty} \widetilde M^3 = 0$.}

\bigskip

Theorem 2 supersedes earlier results by A. Casson, by myself and by many others, from the ninety-nineties, showing that if $\pi_1 M^3$ satisfies any of a long list of geometric condition, like our almost-convexity from Part~I then this implies that $\pi_1^{\infty} \widetilde M^3 = 0$; (see for all this old work, which Theorem 2 supersedes, \cite{3}, \cite{18}, \cite{19}, \cite{20}, but we have only quoted some of the concerned papers here).

\smallskip

At this point, for a locally compact simplicial space $X$, we have seen, here and in Part I, a number of geometric concepts like GSC or WGSC or QSF, which describe its global or asymptotic topological structure. It might be useful for the reader to have a synoptic view of how these various concepts hang together, with their various implications, for a non-compact simplicial complex

\newpage

\centerline{SYNOPTIC VIEW OF ASYMPTOTIC PROPERTIES}

\centerline{for spaces and groups}

$$
\xymatrix{
{\rm GSC} \ar@{=>}[d]^-{{{\rm AFTER \, STABILIZATION} \hfill \atop {\rm \! (see \, the \, stabilization \ lemma \, in \, Part \, I)}}} \ar@{=>}[r] &{\rm WGSC} \ar@{=>}[r] &{\rm QSF} \ar@{<=>}[d]^{{\rm ONLY \, for \, 3\mbox{-}manifolds}} \\
{\mbox{\small Dehn exhaustibility} \atop \mbox{\small(in the simplicial context)}} \ar@{=>}[rr]_-{{\rm ONLY \, for \, 3\mbox{-}manifolds}} &&\pi_1^{\infty} = 0
}
$$

\bigskip

But in this synoptic diagram, it is only QSF which is a clean group-theoretical property too. We also believe that the long quest for $\pi_1^{\infty} \widetilde M^3 = 0$ was a red herring and that the good correct way to phrase ``$\pi_1^{\infty} G = 0$'' is ``$G \in$ QSF''.

\smallskip

The present author has proved the following result, which will be the main topic of the present paper.

\bigskip

\noindent {\bf Theorem 3.} (Po) {\it All finitely presented groups $G$ are QSF.}

\bigskip

The complete proof is contained in the trilogy \cite{21}, \cite{22}, \cite{23} and in this survey we will also refer to these papers as [I], [II], [III], respectively. An early sketchy version of this story is also to be found in \cite{50} but the Trilogy supersedes it. Notice, also, that although QSF makes sense for any locally compact space we, of course, only claim it for finitely presented groups $G$, or equivalently, for the $\widetilde M^3 (G)$ introduced in Part I.

\smallskip

More specifically, our arguments need the finitely presented $G$ and the mantra for the proof of our Theorem 3 is SYMMETRY WITH COMPACT FUNDAMENTAL DOMAIN, and, for the proof in question, both the group action an its related co-compactness, are essential for getting our $\widetilde M^3 (G) \in$ QSF.

\smallskip

Coming now back to our trilogy, what [I] did, was to prove the theorem briefly discussed in Part I (Theorem 4):

\medskip

{\it Any finitely presented group $G$ admits a REPRESENTATION which is locally finite, equivariant, and with uniformly bounded zipping length.}

\medskip

This will be an essential lemma for the main result of [II]. In the present section we will concentrate mainly on [II], trying to explain or give an idea of what is happening there. The core of the proof of Theorem~3, which is also the main result of [II], says that a certain high-dimensional cell-complex $S_u (\widetilde M^3 (G))$ is GSC. This is such that we also have a free action
$$
G \times S_u (\widetilde M^3 (G)) \longrightarrow S_u (\widetilde M^3 (G)) \, ,
$$
but the quotient is non-compact; otherwise we would be already done. There is still the implication
$$
S_u (\widetilde M^3 (G)) \in {\rm GSC} \Longrightarrow G \in {\rm QSF}
$$
to prove. This is the object of part [III] of the trilogy, which we will discuss in the next section of our present survey.

\smallskip

The proof of Theorem 2 is a corollary of the Geometrization of 3-manifolds, proved by Perelman via the Ricci flow, and as already said it also follows from the work of Agol, Wise, and all. The proof of Theorem 3 is totally disjoined and unrelated to that. Of course, also, in the context of Gromov's random groups \cite{45}, \cite{46}, fundamental groups of closed 3-manifolds $\pi_1 M^3$ are very rare events among all the finitely presented groups $G$, with which Theorem 3 deals. And Theorem 3, as such, has not much to do with 3-manifolds.

\smallskip

But then, also, short of invoking Theorem 3 for having Theorem 2 above, as such, the whole big celebrated work of Grisha Perelman, or Agol-Wise, is required. And this means that, Theorem 3 which states a fact valid for {\bf all} finitely presented groups $G$, has to be highly non-trivial, something which contradicts the standardly accepted wisdom.

\smallskip

We believe that the paradox here is only an apparent one and that, similar to the position of rational numbers amongst all the real ones or of quasi-periodic functions among the periodic ones, or of quasicrystals among crystals or of Penrose-type tilings among all tesselations, finitely presented groups are only a small part of a larger category, connected presumably with the non-commutative geometry of Alain Connes, where the paradox should get resolved. But this is only a vague project which we have to leave for later.

\smallskip

The next item should serve as an introduction for the strategy of the proof of Theorem 3. In Part I of this survey we have introduced the notion of {\bf easy} group. With this here is the

\bigskip

\noindent {\bf Lemma 4.} (Otera-Po \cite{40}) {\it Any easy group is QSF.}

\bigskip

We will sketch here the argument for the proof, making use of the 2$^{\rm d}$ version of the definition of easy REPRESENTATIONS as presented in Part I of this survey. But then, combining the Lemma 4 with the fact, proved in Part I that $G$ almost convex implies that $G$ admits an easy REPRESENTATION, we get the following result:

\bigskip

{\it All almost-convex groups, hence the hyperbolic, NIL, a.s.o. are QSF.}

\bigskip

Of course, this immediately follows from Theorem 3, but using Lemma 4 we get it much much cheaper. Keep in mind that the ``easy'' above is not the usual ``easy''.

\smallskip

But, then, we also CONJECTURE that the converse of Lemma 4 is true too, i.e. that all finitely presented groups $G$ are easy. Such a $G$ should always be able to avoid the Whitehead nightmare. Partial results are discussed in Part I of our survey.

\smallskip

So back to Lemma 4 now. And we will assume here that we are given a 2$^{\rm d}$ locally finite REPRESENTATION
$$
X^2 \xrightarrow{ \ f \ } \widetilde M^3 (G) \, , \leqno (1.4)
$$
s.t. both $M^2 (f) \subset X^2 \times X^2$ and $fX^2 \subset \widetilde M^3 (G)$ are closed subsets, inside the respective targets. Nothing like equivariance or bounded zipping length is assumed now. With this, let now
$$
X^2 \equiv X_0^2 \to X_1^2 \to X_2^2 \to \ldots \to fX^2 = X^2 / \Phi (f) = X^2 / \Psi (f) \, , \leqno (1.4.1)
$$
be a zipping strategy. We may assume, without loss of generality that each $X_j^2 \to X_{j+1}^2$ is a compact operation, which homotopically speaking is either a homotopy equivalence or the addition of a 2-cell. More precisely still, we may assume that each of our elementary steps $X_j^2 \to X_{j+1}^2$ is one of the $O(i)$ elementary operations from \cite{27}, \cite{28}, or \cite{40}. The (1.4) has all of its mortal singularities of the undrawable type, and the acyclic $X_j^2 \xrightarrow{ \ O(i) \ } X_{j+1}^2$ creates more such singularities, while in the $O(3)$ operation (which are, homotopically speaking addition of 2-cells) two undrawable singularities meet in head on collision and annihilate each other.

\smallskip

One can show that if any $n \geq 5$ is fixed, then we can thicken each $X_i^2$ in (1.4.1) into a {\bf canonically} attached smooth $n$-manifold $\Theta^n (X_i^2)$ and also change (1.4.1) into a sequence of smooth embeddings, each of which is either a smooth compact dilatation of J.H.C. Whitehead or a smooth addition of a handle of index $\lambda = 2$.

\smallskip

So we write down now the thickened version of the (1.4.1),
$$
\Theta^n (X_0^2) \subset \Theta^n (X_1^2) \subset \Theta^n (X_2^2) \subset \ldots \leqno (1.5)
$$
Since $M^2 (f) \subset X^2 \times X^2$ is closed, and hence also $M_2 (f) \subset X^2$, one can put up together the (1.5) and assemble it into a smooth $n$-manifold with large boundary $\underset{\ell = 0}{\overset{\infty}{\bigcup}} \, \Theta^n (X_{\ell}^2)$. Moreover, since $\Theta^n (X_0^2) = \Theta^n (X^2)$ is GSC and since all the steps in (1.5) preserve GSC, the $\underset{\ell = 0}{\overset{\infty}{\bigcup}} \, \Theta^n (X_{\ell}^2)$ itself is GSC.

\smallskip

Because $fX^2 \subset \widetilde M^3 (G)$ is closed, we also have
$$
\underset{\ell = 0}{\overset{\infty}{\bigcup}} \, \Theta^n (X_{\ell}^2) = \Theta^n (fX^2) \, , \ \mbox{a closed subset of $\Theta^n (\widetilde M^3 (G)$).}
$$
The $2^{\rm d}$ REPRESENTATIONS being essentially surjection, we have finally that
$$
\Theta^n (\widetilde M^3 (G)) = \Theta^n (fX^2) + \{\mbox{handles of index 2 and 3}\} \, ,
$$
making that $\Theta^n (\widetilde M^3 (G))$ is GSC. Also with $n \geq 5$, our construction is {\bf canonical} implying that $G$ acts freely on $\Theta^n (\widetilde M^3 (G))$, with a compact quotient. We can invoke now the implication GSC $\Rightarrow$ QSF and Lemma 1, and get that $G \in {\rm QSF}$. $\Box$

\bigskip

\noindent {\bf A Remark.} -- The condition  $n \geq 5$ is essential for this little argument. If we try to force the approach above to $n=4$ then we meet unsuperable obstacles. This tricky issue is discussed in the introduction to [II], the second part of our trilogy. It is the $n \geq 5$ which makes sure that the construction is canonical, and hence that there is a free action of $G$ when we need it. This ends our REMARK and we move to more adult arguments than Lemma 4.

\smallskip

This also ends the introductory part of this section, and we move now to the first subsection.

\bigskip

\subsection{The 2-dimensional representation theorem}

The reason for moving from the $3^{\rm d}$ REPRESENTATION to a $2^{\rm d}$ one, is that now we have a much more transparent view of the zipping and its limit points. This will be essential for our later geometric realization of the zipping in high dimensions.

\smallskip

So, we will go now to a $2^{\rm d}$ REPRESENTATION THEOREM which pushes to dimension 2 the achievements of the theorem described in Part I of this survey, stating that any $G$ admits $3^{\rm d}$ REPRESENTATIONS which are locally-finite, equivariant and of uniformly bounded zipping length. Now, when we go the dimension two for the REPRESENTATION space $X^2$, the all-important double points set $M_2 (f) \subset X^2$ is in full view.  We cannot assume, a priori, that it is closed, but the next theorem achieves the maximum we can say about it, once it is not closed.

\smallskip

Here comes now the full statement of our $2^{\rm d}$ REPRESENTATION theorem. Afterwards, some necessary comments concerning the proof of the $3^{\rm d}$ REPRESENTATION  itself (Theorems 3, 4 in Part I) will be also offered. And then, some hints of how the $2^{\rm d}$ result should be gotten from the $3^{\rm d}$ one will be given too.

\bigskip

\noindent {\bf Theorem 5.} {\it For any finitely presented group $G$, there is a $2^{\rm d}$ REPRESENTATION
$$
X^2 \xrightarrow{ \ f \ } \widetilde M^3 (G) \, , \leqno (1.6)
$$
with the following features.}

\begin{enumerate}
\item[1)] {\it (First finiteness condition.) The $2^{\rm d}$ REPRESENTATION space $X^2$ is a {\ibf locally finite} $2^{\rm d}$ cell-complex. [But, careful, $fX^2$ is {\ibf not} locally finite.]}
\item[2)] {\it ({\ibf Equivariance}.) There is a free action $G \times X^2 \to X^2$ s.t. $f(gx) = gf(x)$.}
\item[3)] {\it (The second finiteness condition.) We consider now lines $\Lambda \subset X^2$ which are transversal to $M_2 (f) \subset X^2$. Such a transversal will be called {\ibf tight}, if inside $X^2$ we cannot find discs $D^2 \subset X^2$ like in Figure~{\rm 1} below.}

\newpage

$$
\includegraphics[width=45mm]{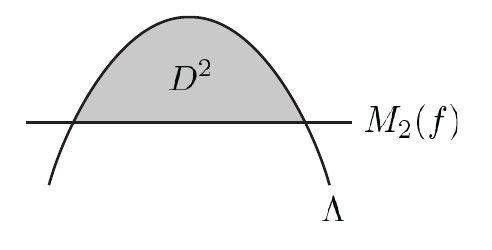}
$$
\label{fig1}
\centerline {\bf Figure 1.} 

\medskip

\centerline{\it A transversal $\Lambda$ to $M_2 (f)$, inside $X^2$, which is {\ibf not tight}.}

\bigskip

{\it With this, our present condition concerns the accumulation points of $\Lambda \, \cap \, M_2 (f)$ and it says the following. For any compact tight transversal $\Lambda$ we have}
$$
{\rm Card} \, (\lim (\Lambda \cap M_2 (f)) < \infty \, .
$$
\item[4)] {\it The following closed subset
$$
{\rm LIM} \, M_2 (f) \equiv \underset{\overbrace{\mbox{\footnotesize $\Lambda \in \{$tight transversals$\}$}}}{\bigcup} \lim (\Lambda \cap M_2 (f)) \subset X^2 \, , \leqno (1.6.1)
$$
which is the only place where the double points $M_2 (f) \subset X^2$ can accumulate, is a locally finite graph.

Moreover the following set
$$
f( {\rm LIM} \, M_2 (f)) \subset fX^2
$$
is {\ibf closed} too. (Notice that, by now, this is not a consequence of general principles, but something which needs its special proof.)}
\item[5)] {\it We also have
$$
\underset{\overbrace{\mbox{\footnotesize $\Lambda \in \{$tight transversals$\}$}}}{\bigcup} \Lambda \cap M_2 (f) = M_2 (f) \, .
$$
}
\item[6)] {\it ({\ibf Uniformly bounded zipping length}.) For any $(x,y) \in M^2(f) \subset X^2 \times X^2$ we consider the {\ibf zipping paths}
$$
\lambda (x,y) \subset \widehat M^2 (f) = M^2 (f) \cup {\rm Sing} \, (f)
$$
(with ${\rm Sing} \, (f) = {\rm Diag} \, ({\rm Sing} \, (f))$). Such zipping paths do have to be there for any $(x,y) \in M^2 (f)$, since {\rm (1.6)} is a REPRESENTATION. There exist a uniform bound $K > 0$ such that, when $\lambda$ runs over all the zipping paths of all double points $(x,y)$, then $\underset{\lambda}{\rm inf}$ length $\lambda (x,y) < K$.}
\end{enumerate}

\bigskip

\noindent {\bf Comments.} -- The ${\rm LIM} \, M_2 (f)$ has to be compared with foliations $({\mathcal F})$ or laminations $({\mathcal L})$, in $2^{\rm d}$, with $1^{\rm d}$ leaves. Like for ${\mathcal L}$, our ${\rm LIM} \, M_2 (f)$ has a holonomy, which may be non-trivial. Figure 2 shows the local $R^2$-models for ${\mathcal F}$, ${\mathcal L}$ and ${\rm LIM} \, M_2 (f)$, and on the right column marks their respective transverse structures.

$$
\includegraphics[width=100mm]{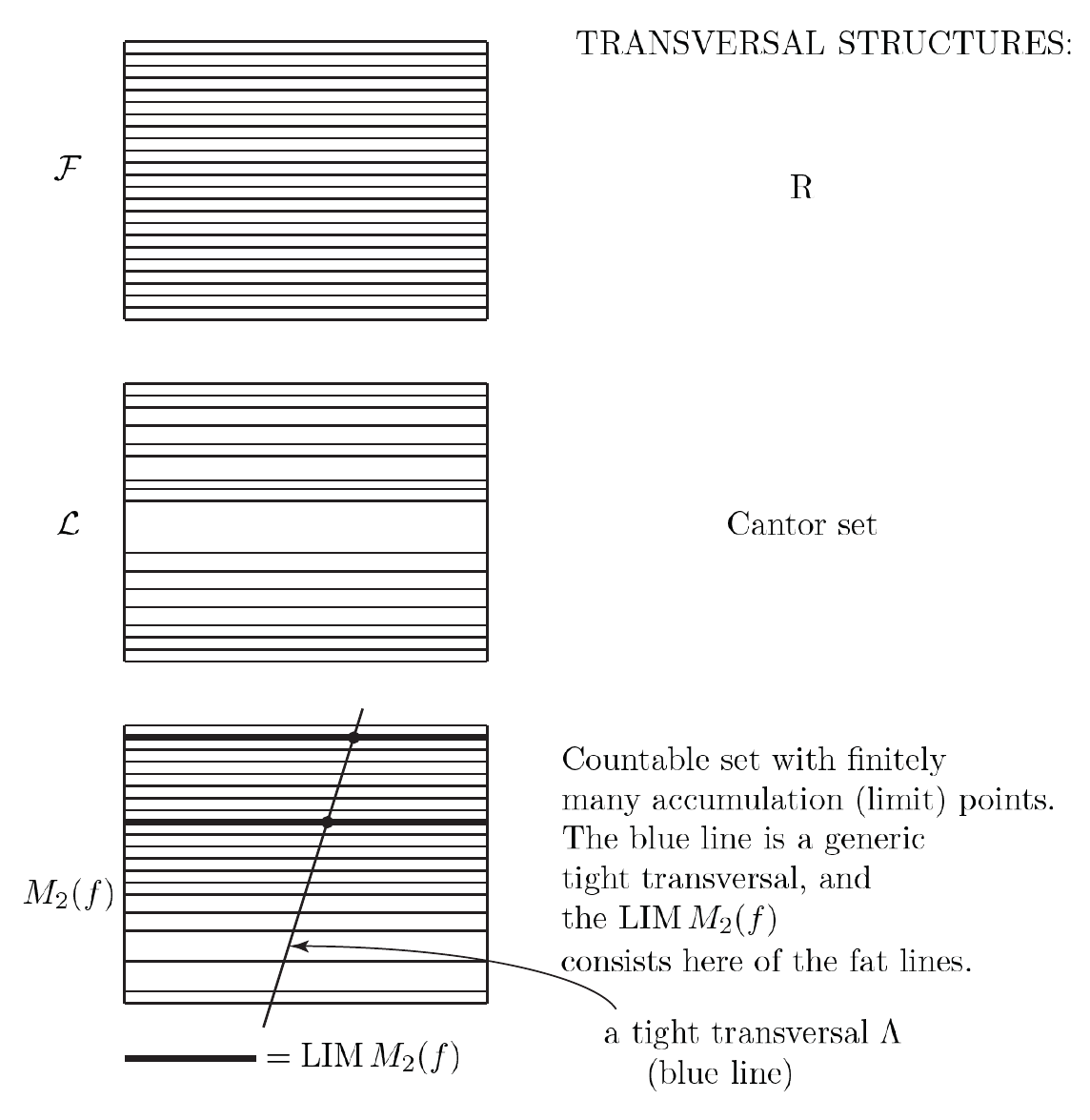}
$$
\label{fig2}
\centerline {\bf Figure 2.} 

\bigskip

\noindent Notice, finally, that the length $\lambda (x,y)$ is well-defined up to a quasi-isometry, making our condition $\underset{\lambda}{\rm inf}$ length $\lambda (x,y) < K$ well-defined. End of Comments.

\bigskip

Figure 2 displays the simplest local model for $(M_2 (f) , {\rm LIM} \, M_2 (f))$, and a more complete model is provided by Figure 3. But even in Figure 3 we are still inside a local chart $R^2 \subset X^2$, far from the non-manifold points in $X^2$ and from its mortal singularities. At this point, a few things will also be said about how the local finiteness of $X^3$ is achieved for the $3^{\rm d}$ REPRESENTATION $X^3 \xrightarrow{ \ f \ } \widetilde M^3 (G)$ which has all those nice features of local finiteness, equivariance and uniformly bounded zipping length.

$$
\includegraphics[width=120mm]{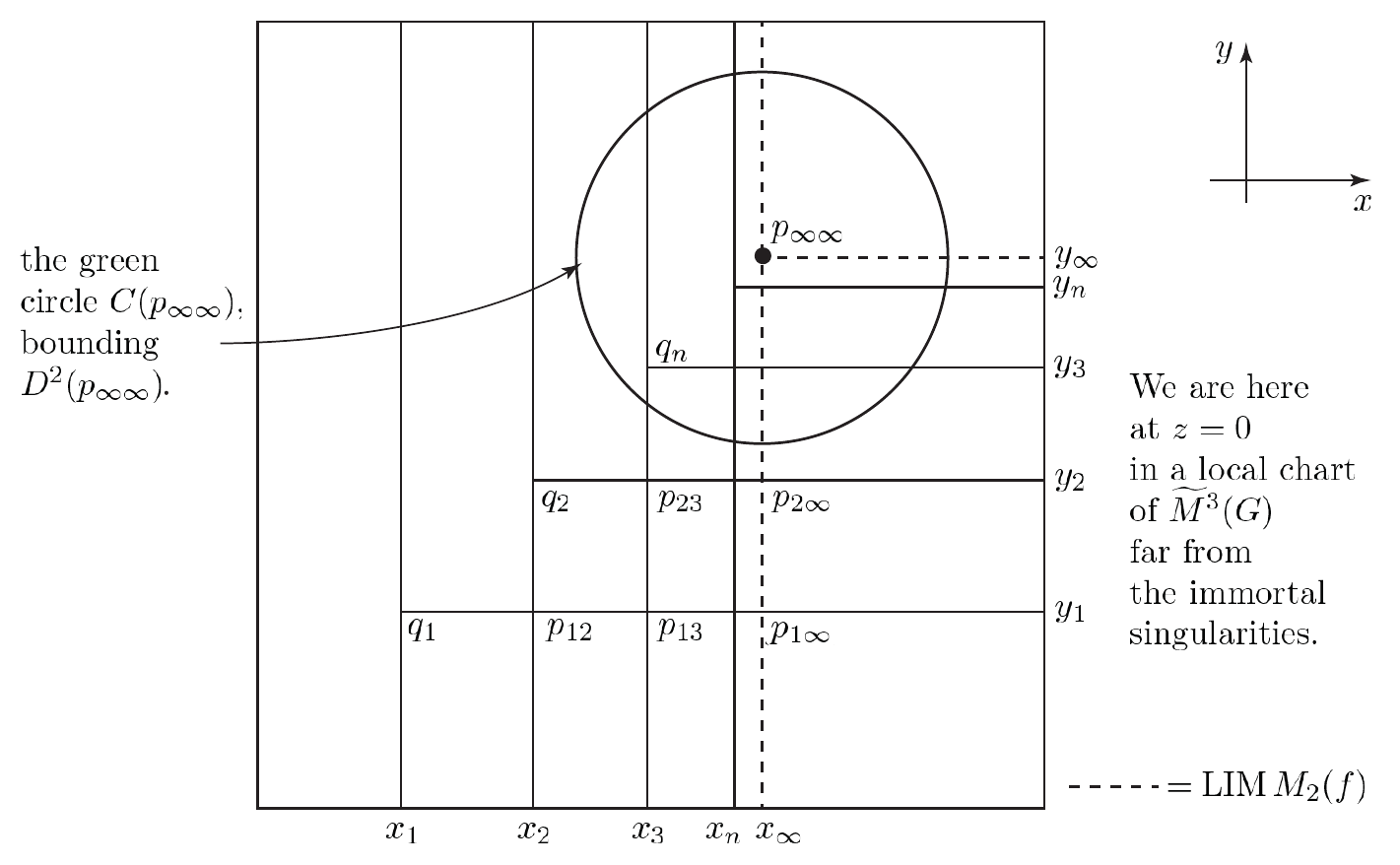}
$$
\label{fig3}
\centerline {\bf Figure 3.} 
\begin{quote}
A local model of $M_2 (f) \subset X^2$ and its accumulations. Think here of the background square as being BLACK, actually, with the notations to be introduced later, it is a $W_{(\infty)}$ (BLACK) piece of $X^2$, each vertical plane $(x = x_i , y , z)$ a BLUE sheet of $X^2$ and each plane $(x,y = y_j , z)$ a RED sheet. The exact meaning of these colours will be soon explained. A vertical line along a $q_i$-point is a line along which $X^2$ is not a 2-manifold and each $p_{ij} = (x = x_i , y = y_j)$ is the trace of a line in $M_2 (f)$, transversal to the plane of the figure, in $\widetilde M^3 (G)$. This also means that each $p_{ij}$ is a triple point, living in $M_3 (f) \subset X^2$. The planes $(x=x_{\infty} , y , z)$, $(x,y = y_{\infty} , z)$ are limiting positions of BLUE or RED sheets in $X^2$, respectively. Their dotted trace on the BLACK square is in ${\rm LIM} \, M_2 (f)$.

For reason of simplicity in this very localized figure we assume to be here far from any kind of singularities, mortal $({\rm Sing} \, (f))$ or immortal $({\rm Sing} \, \widetilde M^2 (G))$. As one might guess, the points denoted $p_{\infty\infty}$ play a special role; they are source of infinite headaches, which will have to be taken care of.
\end{quote}

\bigskip

We consider now the context of the naive, Kindergarten theory of universal covering spaces and the pristine $3^{\rm d}$ GOOD REPRESENTATIONS, which were introduced at that occasion.

\smallskip

We are there in the presence of two $3^{\rm d}$ singular handlebodies
$$
\widetilde M^3 (G) = \underset{\overbrace{\mbox{\footnotesize $0 \leq \lambda \leq 3 , i$}}}{\sum} h_i^{\lambda} \, , \quad X^3 = \underset{\overbrace{\mbox{\footnotesize $\lambda , i , \alpha$}}}{\bigcup} h_i^{\lambda} (\alpha) \, , \leqno (1.7)
$$
where $\alpha$ is an index belonging to some countable family and, where for each given $i,\lambda$ we have an isomorphism $f (h_i^{\lambda} (\alpha)) = h_i^{\lambda}$. 

\smallskip

We move now to the more interesting context of the existence, for each $G$, of a $3^{\rm d}$ REPRESENTATION which is locally finite, equivariant and with uniformly bounded zipping length. This was a theorem stated in Part I of this survey. So, we have now a REPRESENTATION $X^3 \xrightarrow{ \ f \ } \widetilde M^3 (G)$ like before. Here $\widetilde M^3 (G)$ is like above but $X^3$ is more sophisticated. It is infinitely larger than before, a price to be paid for $X^3$ to be locally finite and with a free action $G \times X^3 \to X^3$. Its building blocks are the so-called {\bf bicollared handles}, which we will introduce now, but see also \cite{21}, \cite{29}.

\bigskip

Topologically speaking, a bicollared handle of dimension $n$ and index $\lambda$ is a copy of $B^{\lambda} \times {\rm int} \, B^{n-\lambda}$, with the lateral surface $\delta H^{\lambda} = B^{\lambda} \times \partial \, B^{n-\lambda}$ living at infinity, and with additional structures to be explained. Of course, we are only interested in the case $n=3$ but we find the more general notation easier to follow. Here is a bicollared structure, unrolled.

\smallskip

We start with $R^n = R^{\lambda} \times R^{n-\lambda}$ and with two infinite sequences of concentric balls

\smallskip

$R^{\lambda} \supset B_1^{\lambda} \supset B_2^{\lambda} \supset \ldots \supset B_m^{\lambda} \supset \ldots$, bounded from below, and

\smallskip

$B_1^{n-\lambda} \subset B_2^{n-\lambda} \subset \ldots \subset B_m^{n-\lambda} \subset \ldots \subset R^{n-\lambda}$, such that $\underset{m = \infty}{\rm lim} B_m^{n-\lambda} = R^{n-\lambda}$.

\smallskip

With this, for every $i=1,2,\ldots$ we consider the standard $\lambda$-handle $H_i^{\lambda} = B_i^{\lambda} \times B_i^{n-\lambda}$ and, by definition, the
$$
H^{\lambda} = \bigcup_{i=1}^{\infty} H_i^{\lambda} \qquad \mbox{(see Figure 4)}, \leqno (1.8)
$$
is a bicollared handle of index $\lambda$. Figure 4 suggests for our $H^{\lambda}$ two {\bf collars} as well as the attaching zone
$$
\partial H^{\lambda} \cong \partial B^{\lambda} \times {\rm int} \, B^{n-\lambda}
$$
which is suggested in fat lines. By ``collar'' we mean something of the form $\partial Y^n \times [1,\infty)$ together with the filtration
$$
\partial Y^n \times [1,\infty) \supset \partial Y^n \times [2,\infty) \supset \partial Y^n \times [3, \infty) \ldots
$$
and with the {\ibf k}-{\bf levels} $\partial Y^n \times [k,k+1]$.

$$
\includegraphics[width=110mm]{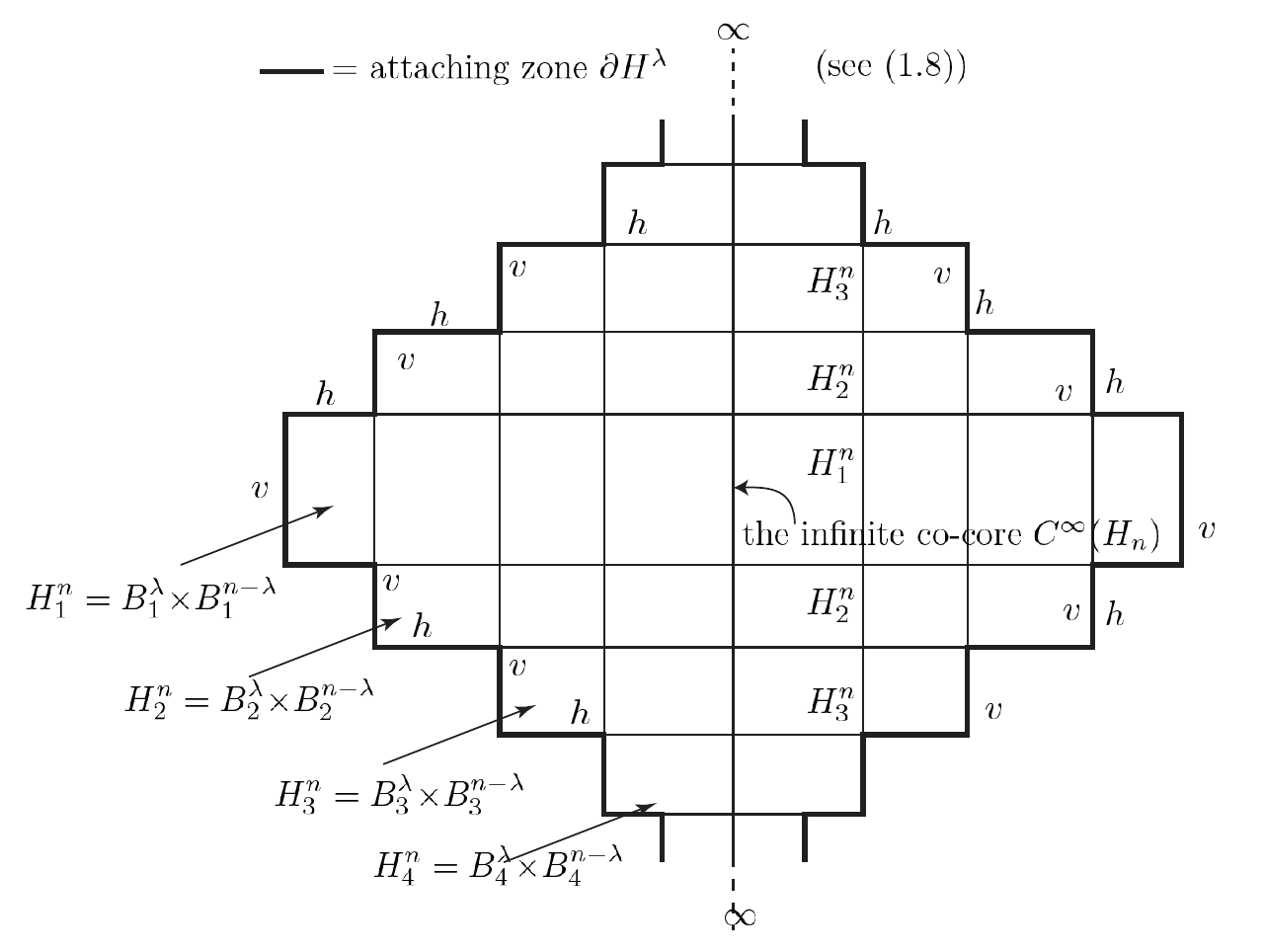}
$$
\label{fig4}
\centerline {\bf Figure 4.} 
\begin{quote}
The bicollared handle $H^n (\lambda) = \bigcup_1^{\infty} B_i^{\lambda} \times B_i^{n-\lambda}$. The attaching zones of the $H_i^n$'s stay far from the infinite co-core $C^{\infty} (H^n) = \bigcup_m$ co-core $(H_m^n)$. Let us say that they have limiting positions parallel to $C^{\infty} (H^n)$ and disjoined from it. Here $h, v$ stand for ``horizontal'' and ``vertical'', respectively.
\end{quote}

\bigskip

The bicollared handle has an ingoing collar, starting at the attaching zone $\partial H^n$ and an outgoing collar, going towards the missing lateral surface; these collars are not disjoined, of course.

\smallskip

When, for $\lambda_1 > \lambda_2$, $H^{\lambda_1}$ is attached to $H^{\lambda_2}$, we use the ingoing collar of $H^{\lambda_1}$ and the outgoing one of $H^{\lambda_2}$, starting at some {\bf chosen level} $k \in Z_+$. In order to achieve equivariance we will need to attach infinitely many $H_1^{\lambda_1} , H_2^{\lambda_1} , \ldots$ to a same $H^{\lambda_2}$. When that will happen we will use levels $k_1, k_2 \ldots$ with $\underset{n=\infty}{\lim} \, k_{n} = \infty$, and this will yield our local finiteness of $X^3$. Our $X^3$ is $X^3 = \underset{i , \gamma , \lambda \leq 3}{\bigcup} \, H_i^{\lambda} (\gamma)$, with $\{i\}$ the same family of indices as in (1.7), and with $\{\gamma\}$ another infinite family of indices, extending $\{\alpha\}$.

\smallskip

In terms of (1.7), each bicollared $H^{\lambda}$ is a souped up version of some $h_i^{\lambda} (\alpha)$ in (1.7). But careful, the $X^3$ in (1.7) is the pristine, naive one, now we talk about the real life $X^3$, with all those fancier features like local finiteness, equivariance, a.s.o.

\smallskip

For the real-life map $f$ considered now, the $f H_i^{\lambda} (\gamma)$ (now for a much larger family of indices $\gamma$ than the pristine $\alpha$'s) occupies, in $\widetilde M^3 (G)$ {\it roughly} the position $h_i^{\lambda}$, but in such a way that $f \mid H_i^{\lambda} (\gamma)$ extends continuously to the completion

\bigskip

\noindent (1.8.1) \quad $\widehat H_i^{\lambda} (\gamma) = H_i^{\lambda} (\gamma) \cup \{$the mythical $\delta H_i^{\lambda} (\gamma)$, which lives at the infinity of $H_i^{\lambda} (\gamma)\} \underset{\rm TOP}{=} B^{\lambda} \times B^{n - \lambda}$,

\bigskip

\noindent in such a way that we should have {\bf strict} equality, for each $i,\gamma$
$$
f (\delta H_i^{\lambda} (\gamma)) = \delta h_i^{\lambda} \, .
$$

In order to get our desired features of local-finiteness and equivariance for our real-life $3^{\rm d}$ REPRESENTATION space $X^3$, it is necessary that when one moves from the families $\{h_i^{\lambda} (\alpha) , (\lambda , i) \}$, in (1.7) to the corresponding $\{ H_i^{\lambda} (\gamma) , (\lambda , i)\}$ of bicollared handles we vastly increase the family of indices $\{\alpha\}$, into $\{\gamma\} \supsetneqq \{\alpha\}$.

\smallskip

We will not say more here about how to get that $3^{\rm d}$ REPRESENTATION $X^3 \xrightarrow{ \ f \ } \widetilde M^3 (G)$ with its local finiteness, equivariance a.s.o. So we move from it, which we assume given, to the context of our Theorem~5. But, for pedagogical reasons, we will start with the simpler context of the naive pristine $X^3 \xrightarrow{ \ f \ } \widetilde M^3 (G)$ from the Kindergarten theory of universal covering spaces, without any more sophisticated adult features (and see here Part I of this survey), and we will show how to turn {\bf it} into a naive $2^{\rm d}$ REPRESENTATION  $X^2 \xrightarrow{ \ f \ } \widetilde M^3 (G)$, forgetting temporarily about the items 1) to 6) in our Theorem 5.

\smallskip

One starts by foliating, compatibly, the $3^{\rm d}$ handles of $\widetilde M^3 (G)$ and $X^3$. For each handle $h^{\lambda}$, where $\lambda = 0,1,2$, and we do not need to worry much about $\lambda = 3$, as it turns out, three {\bf not everywhere defined} foliations will be considered. Call them ${\mathcal F}$ (COLOUR) and these colours will be BLUE, RED, BLACK which will be ``{\bf natural}'' (a word to be explained) for $\lambda = 0$, $\lambda = 1$, $\lambda = 2$ respectively. Our foliations always have leaves of dimension two. For the $h^{\lambda}$ and for its natural colour ${\mathcal F}={\mathcal F}$ (COLOUR, natural for $\lambda$), the leaves are all isomorphic copies of the lateral surface $B^{\lambda} \times \partial B^{n-\lambda}$, parallel to it and coming closer and closer to the core $B^{\lambda} \subset h^{\lambda}$ which is not part of any leaf of the ${\mathcal F}$ (natural $\lambda$-COLOUR). For each individual handle $h^{\lambda} (\alpha) \subset X^3$, each leaf of ${\mathcal F}$ (COLOUR), when restricted to $h_i^{\lambda} (\alpha)$ is called a {\bf compact wall}, and for given $\lambda$ and for the corresponding NATURAL COLOUR these are, respectively, copies of $S^2$ (BLUE) (for $\lambda = 0$), $S^1 \times [0,1]$ (RED) (for $\lambda = 1$) and $D^2$ (BLACK) (for $\lambda = 2$).

\smallskip

For any $h^{\lambda}$ all the tree ${\mathcal F}$ (COLOUR) make sense, one colour, exactly, being natural, the other two unnatural. So, we have here, inside each $h^{\lambda}$, a triply orthogonal system of compact walls, not everywhere defined, of course. With this, the colours BLACK, RED, BLUE which were already mentioned in the context of Figure 3 should already make sense. The general idea in going from $X^3$ to $X^2$ is to build $X^2$ as an infinite union of compact walls, which should make up a very dense 2-skeleton of $X^3$.

\smallskip

This idea of ${\mathcal F}$ (COLOUR) extends to the bicollared handles and to the sophisticated, adult $X^3$, which we use now. But now we need to take care too of the items 1) to 6) in Theorem 5 too. And to the various conditions expressed in the theorem in question, we also want to add now the following one too

\newpage

\noindent (1.9) \quad ``In principle'' (and the quotation marks will be clarified later on), the mortal singularities of $X^2$ are all of the undrawable type as described in Part I of this survey (see (1.10.1) in Part I) and, more importantly
$$
\{\mbox{the set of undrawable mortal singularities of (1.6)}\} \subset X^2
$$
is {\bf discrete}, i.e. without accumulation at finite distance. End of (1.9).

\bigskip

The (1.9), as well as other constraints, force us to put into the $X^2$, $2^{\rm d}$ REPRESENTATION space occurring in (1.6), more than just the compact walls. We will also need to add {\bf non-compact security walls} $W_{\infty} \subset X^2$, all of them of BLACK colour. And now, very importantly, the BLACK square in Figure 3 could well live inside such a $W_{\infty}$, or inside a compact, also BLACK wall. We will denote by $W_{(\infty)}$ (BLACK) a black wall which is either a compact  $W$ (BLACK) or a non-compact security wall.

\smallskip

All this being said, here is how the ${\rm LIM} \, M_2 (f)$ (1.6.1) in Theorem 5 occurs. We have

\bigskip

\noindent (1.10) \quad $\Sigma_1 (\infty) \equiv \{$the union of the limit positions of the {\bf compact} walls $fW\} \subset \widetilde M^3 (G)$,

\bigskip

\noindent and with this, we also have now
$$
{\rm LIM} \, M_2 (f) = f^{-1} (fX^2 \pitchfork \Sigma_1 (\infty)) \subset X^2 \, , \leqno (1.11.1)
$$
$$
f \, {\rm LIM} \, M_2 (f) = fX^2 \pitchfork \Sigma_1 (\infty) \subset fX^2 \, . \leqno (1.11.2)
$$
It turns out that the BLACK limiting positions do not contribute to the $\Sigma_1 (\infty)$ in (1.11.1) and (1.11.2) and not the $W_{(\infty)}$ (BLACK) either.

\smallskip

Figure 3 actually presents the typical instance of (1.11.1) and all the pieces of $\Sigma_1 (\infty)$ involved there are BLUE or RED, cutting through a BLACK wall (which can be compact or a security wall).

\smallskip

There is more to be said concerning Figure 3. We see there also infinitely many zig-zag polygonal lines going through successive $q$'s and $p$'s and some of them are not double points but branching loci of $X^2$ created by the attachments of the bicollared handles. And with this, the presence of the $p_{\infty \infty}$ would make our $X^2$ not be locally finite. Of course, also, a very violent violation of the local finiteness of $fX^2$ at $p_{\alpha \infty}$ occurs now. These things force us to go to a more sophisticated way of conceiving the $X^2$ (and also $fX^2$), and so we will go from the $X^2$ as presented so far, call it the naive $X^2$ to the following object which is now locally finite.

\bigskip

\noindent (1.12) \quad When we go the the real life $X^2$ we delete $p_{\infty\infty}$ (from the naive $X^2$) and we compensate this by adding the disc $D^2 (p_{\infty\infty})$ (Figure 3) along the circle $C(p_{\infty\infty})$. This $C(p_{\infty\infty})$ becomes now a new singularity which is a folding line. Of course, this is now no longer one of our undrawable singularities. This should explain the quotation marks in the context of (1.9). The fold-singularities $C(p_{\infty\infty}) = S^1$ do not accumulate at finite distance and they do not create any serious difficulties for us.

\bigskip

So far, we have totally ignored the immortal singularities of $\widetilde M^3 (G)$. When one looks in more detail at what happens there, at $X^2$-level, there are more $p_{\infty\infty}$-type nasty points, actually quite nastier, occurring in connection with the immortal singularities of $\widetilde M^3 (G)$. These nastier singularities $p_{\infty\infty} (S) \in X^2$ (not to be mixed up the $p_{\infty\infty}$ (proper), from Figure 3 and (1.12)), are there iff $M^3 (G)$ has immortal singularities, i.e. if and only if $G$ is not the $\pi_1$ of a compact smooth 3-manifold $M^3$. They will have to be deleted too, like in (1.12) but in order to keep this survey at reasonable length, we will not discuss them very much right now.

\smallskip

Let us go back now to our $M_2 (f)$ and to its associated ${\rm LIM} \, M_2 (f)$, as displayed in Figure 3.

\smallskip

It is very important that in the group-theoretical set-up of Theorem 5, the accumulation pattern for ${\rm LIM} \, M_2 (f)$ is ruled by the card $(\lim (\Lambda \cap M_2 (f)) < \infty$ from 3) in our theorem. But the ${\rm LIM} \, M_2 (f)$ game can be played in other contexts too, and when one looks at its accumulation pattern in the context of the wild Whitehead manifolds ${\rm Wh}^3$, then the Julia sets of the complex dynamics pop up (see \cite{30}), i.e. chaos appears.

\smallskip

Notice, also that we have the double implication
$$
{\rm LIM} \, M_2 (f) \ne \emptyset \Longleftrightarrow M_2 (f) \subset X^2 \ \mbox{is NOT closed.}
$$

So, once we are not in the case of an easy REPRESENTATION, the ${\rm LIM} \, M_2 (f)$ has to be there. And our condition
$$
{\rm card} \, (\lim (\Lambda \cap M_2 (f)) < \infty
$$
is then the most benign of the possible evils which might occur.

\smallskip

The next item is an important complement to Theorem 5.

\bigskip

\noindent (1.13) \quad For the equivariant REPRESENTATION (1.6), we can find an {\bf equivariant zipping}
$$
\qquad\qquad\qquad X^2 = X_0^2 \to X_1^2 \to X_2^2 \to \ldots \to fX^2 = X^2 / \Phi (f) \, ,
$$
\vglue-6mm
$$
\mbox{\hglue 7mm} {\mid\mbox{\hglue-1mm}}_{\!\!-\!\!-\!\!-\!\!-\!\!-\!\!-\!\!-\!\!-\!\!-\!\!-\!\!-\!\!-\!\!-\!\!-\!\!-\!\!-\!\!-\!\!-\!\!-\!\!-\!\!-\!\!-\!\!-\!\!-\!\!-\!\!-\!\!-\!\!-\!\!-\!\!-\!\!-\!\!-\!\!-\!\!-\!\!-\!\!-\!\!-\!\!-\!\!-\!\!-\!\!-\!\!-\!\!-\!\!-\!\!-\!\!-}{\mbox{\hglue -2mm}\uparrow}
$$
\vglue-5mm
$$
\qquad f
$$
giving rise to the commutative diagram
$$
\xymatrix{
X^2 \ar[d] \ar[rrr]^-f_-{\mbox{\footnotesize equivariant zipping,} \atop \mbox{ \footnotesize upstairs}} &&&\widetilde M^3 (G) \, . \ar[d] \\
X^2 / G \ar[rrr]^-{f/G}_-{\mbox{\footnotesize equivariant zipping,} \atop \mbox{ \footnotesize downstairs}} &&&M^3 (G)
}
$$
Here, the upper line is REPRESENTATION (of $G$), but certainly not the lower one since $\pi_1 (X^2/G) \ne 0$ (hence $X^2/G$ cannot be GSC). Nor are the partial
$$
X^2_{n < \infty} \underset{ \ f \mid X_n^2}{-\!\!\!-\!\!\!-\!\!\!-\!\!\!-\!\!\!\longrightarrow} \widetilde M^3 (G) \, ,
$$
which are certainly not GSC. End of (1.13).

\bigskip

From now on, in this section and the next, only equivariant zipping will be considered. This ends our discussion of Theorem 5, and after the $2^{\rm d}$ interlude we go back to $3^{\rm d}$, and we open another subsection.

\subsection{Where we move back to dimension three}

With things like the procedure (1.12) our $X^2$ is locally finite. It is certainly GSC too and we will consider for it a canonical $3^{\rm d}$ thickening $\Theta^3 (X^2)$, a $3^{\rm d}$ {\bf singular manifold} with very large boundary and also a smooth high-dimensional thickening of $\Theta^3 (X^2)$ denoted by
$$
\Theta^{N+4} (X^2) = \Theta^4 (\Theta^3 (X^2)) \times B^N \, , \ \mbox{{\bf not a manifold} either.} \leqno (1.14)
$$
Here ``$\Theta$'' stands for ``thickening'' and the general idea behind a formula like (1.14) (and more like this will come), is the following, and further details concerning $\Theta^3 (X^2)$ will follow now.

\smallskip

We want now a singular ``$3^{\rm d}$ neighbourhood'' $\Theta^3 (X^2)$ of $X^2$, which, on the other hand will help make good sense of things like (1.12) but, more seriously, which will be a first step for the GEOMETRIC REALIZATION OF THE ZIPPING, which will happen in a very high dimension $N+4$.

\smallskip

Some words of caution are here in order. To begin with, $\Theta^3 (X^2)$ is NOT just another version of $X^3$ and its virtue is that it offers a better grip on the zipping. And then, our way of using those high dimensions, has nothing to do with the well-known high dimensions in differential topology (see Smale's $h$-cobordism theorem and all that). Also, this singular, locally finite object, is {\bf not} a regular neighbourhood of $X^2$, whether we mean that naive $X^2$ or its more sophisticated version. This is the reason for the quotation marks above.

\smallskip

The $\Theta^3 (X^2)$ is a 3-manifold with singularities of the undrawable type and also simple-minded fold like singularities coming from the $C(p_{\infty\infty})$. Outside these latter ones, which we want to leave them {\bf singular}, forever, our $\Theta^3 (X^2)$ can be smoothened (and this concerns now only the undrawable singularities by going to dimension four, i.e. there is a smooth $4^{\rm d}$ thickening $\Theta^4 (\Theta^3 (X^2))$. The problem is that this is certainly not unique, it depends on certain choices; and see here the desingularizations ${\mathcal R}$ (2.6), which are discussed in the next section. And see here \cite{27}, \cite{28} too. It turns out that by multiplying with $B^N$ we wash out the difference between these choices, coming with ${\mathcal R}$, and $\Theta^4 (\Theta^3 (X^2)) \times B^N$ (1.14) becomes {\bf canonical}. And, because it is canonical, the free action $G \times X^2 \to X^2$ extends now to a free action
$$
G \times \Theta^{N+4} (X^2) \to \Theta^{N+4} (X^2) \, .
$$

Coming back to (1.14), our $\Theta^4 (\Theta^3 (X^2))$ and $\Theta^{N+4} (X^2)$ will continue to be singular at the $C(p_{\infty\infty})$, because in the later part of our arguments we will want to be able to dispose freely of the 2-handles $D^2 (p_{\infty\infty})$ which compensate for the gone $p_{\infty\infty}$, and in particular of their attaching curves.

\smallskip

In a more serious vein this kind of little game will also function for the $(N+4)$-dimensional objects $S_u (\widetilde M^3 (G))$ and $S_b (\widetilde M^3 (G))$ to be encountered later on.

\smallskip

There is in all this an important distinction between the ``{\bf high dimensions}'', $N+4$ in our present set-up, similar to the ones of Smale and Stallings and the ``{\bf supplementary dimensions}'', $N$ in addition to the four of $\Theta^4 (\Theta^3 (X^2))$. It is there that our main action will take place.

\smallskip

Finally, $\Theta^{N+4} (X^2)$ is certainly GSC, just like $\Theta^3 (X^2)$.

\smallskip

So far so good, but even with (1.12), the $fX^2$ is still not locally finite, and this comes from ${\rm LIM} \, M^2 (f) \supset \{p_{\infty\infty}\}$, see here the Figure 3, but keep in mind too that there are even worst $p_{\infty\infty}$'s occurring in connection with the immortal singularities of $\widetilde M^3 (G)$, the $p_{\infty\infty} (S)$, as opposed to $p_{\infty\infty}$ (proper) from Figure 3; we will have to come back to them.

\smallskip

One should notice that when one goes from $\{$the sophisticated $X^2\} \equiv \{$naive $X^2\} - \Sigma \, p_{\infty\infty} + \underset{p_{\infty\infty}}{\Sigma} \{$the compensating disc $D^2 (p_{\infty\infty})$, of center $p_{\infty\infty}\}$ to $fX^2$, then $p_{\infty\infty}$, now center of $D^2 (p_{\infty\infty})$, is still there as a point of non-local finiteness of $fX^2$. We will nevertheless introduce a {\bf locally finite} $3^{\rm d}$ object $\Theta^3 (fX^2)$. This will be a $3^{\rm d}$ singular manifold analogous to $\Theta^3 (X^2)$ and which, for us, will be the correct $3^{\rm d}$ thickening of the $fX^2$; but there will be NO inclusion $fX^2 \subset \Theta^3 (fX^2)$. Here is how this goes. Having one additional dimension available, the lack of local finiteness of $fX^2$ coming from
$$
fX^2 \supset f \, {\rm LIM} \, M_2 (f) \supset \{p_{\infty\infty}\}
$$
will be taken care of by {\bf deletions}  which will be either sending certain lines of the boundary to infinity (see the fat red points in the Figure 5) or punching holes which we compensate with 2-handles, like in the Figure~6. These 2-handles are singular, and they stay singular at all the levels in formula (1.14). That is why the objects there stay not smooth. When one is far from $p_{\infty\infty}$ and from the immortal singularities of $\widetilde M^3 (G)$, we have the following kind of local model for $X^2 \xrightarrow{ \ f \ } \widetilde M^3 (G)$; and in order to simplify our exposition, we will stay far from the $M_3 (f)$ too. There is a local chart $U = R^3 = (x,y,z) \subset \widetilde M^3 (G)$ inside which live the $\infty + 1$ many planes $W = (z=0)$ and $V_n = (x=x_n)$, where $x_1 < x_2 < \ldots$, $\underset{n = \infty}{\lim} \, x_n = x_{\infty}$. The local model for the $2^{\rm d}$ REPRESENTATION $X^2 \xrightarrow{ \ f \ } \widetilde M^3 (G)$, is $f^{-1} U = W \cup \underset{1}{\overset{\infty}{\Sigma}} \, V_n \subset X^2$, with the obvious maps $f \mid W$, $f \mid V_n$. Here $(x=x_{\infty} , z=0) \subset W$ is $f({\rm LIM} \, M_2 (f)) \cap U$ and, in this very simple situation, it is not necessary to distinguish between ${\rm LIM} \, M_2 (f)$ and $f({\rm LIM} \, M_2 (f))$.

\smallskip

For $\Theta^3 (fX^2)$, the corresponding local model (which is non-singular in this case) is
$$
\{\mbox{local, non-singular, locally finite piece of} \ \Theta^3 (fX^2)\} = \left[ W \times [-\varepsilon \leq z \leq \varepsilon] - {\rm LIM} \, M_2(f) \times (z = \pm \, \varepsilon)\right] \, \cup
$$
$$
\cup \ \sum_{n=1}^{\infty} V_n \times (x_n - \varepsilon_n \leq x \leq x_n + \varepsilon_n) \leqno (1.15)
$$
where the $\varepsilon_n$'s are converging very fast to zero. The (1.15) is suggested in the Figure 5, which is a cross-section $y = {\rm const}$ through $U \cap \Theta^3 (fX^2)$.
$$
\includegraphics[width=95mm]{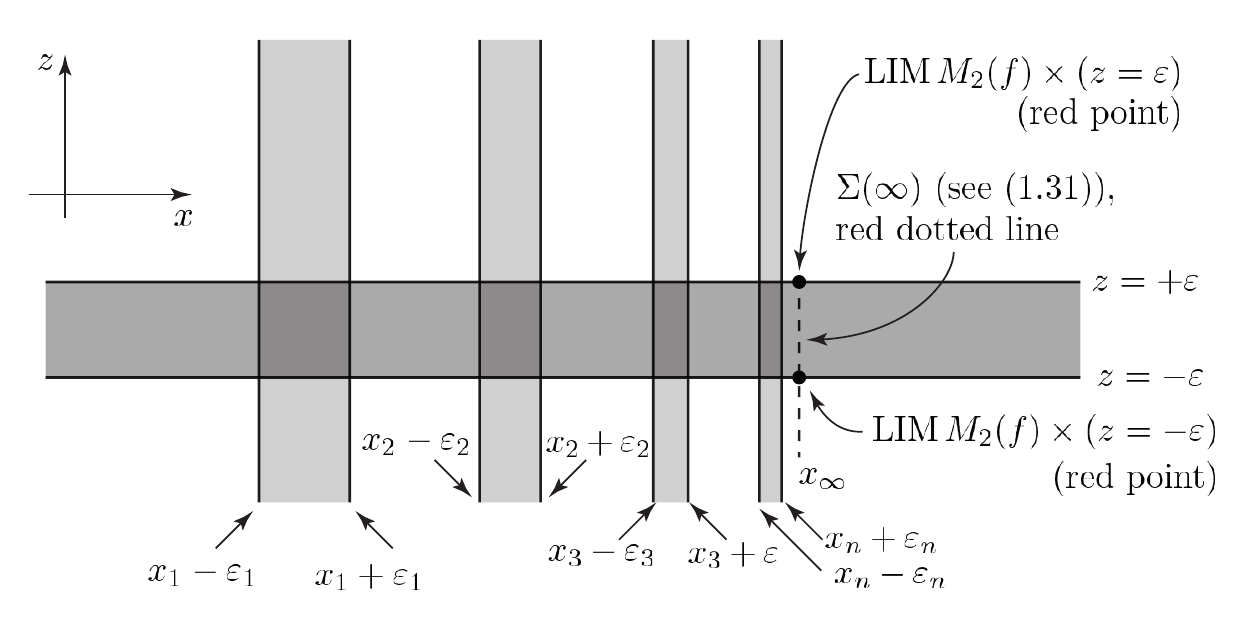}
$$
\label{fig5}
\centerline {\bf Figure 5.} 
\begin{quote}
A section through $\Theta^3 (fX^2)$. The fat RED lines ${\rm LIM} \, M_2(f) \times (z = \pm \, \varepsilon)$ are deleted, so as to stay locally finite. Here, the shaded region $- \varepsilon \leq z \leq \varepsilon$ is a thickened $W$ (BLACK) or $W_{\infty}$ (BLACK), while all the vertical $\{x_n - \varepsilon_n \leq x \leq x_n + \varepsilon_n \}$'s get thinner and thinner as $n \to \infty$. They are all $W$ (BLUE) or $W$ (RED), thickened.
\end{quote}

\smallskip

All this was far from the $p_{\infty\infty}$'s and when we come close to them we proceed like it is suggested in the Figure 6.

$$
\includegraphics[width=120mm]{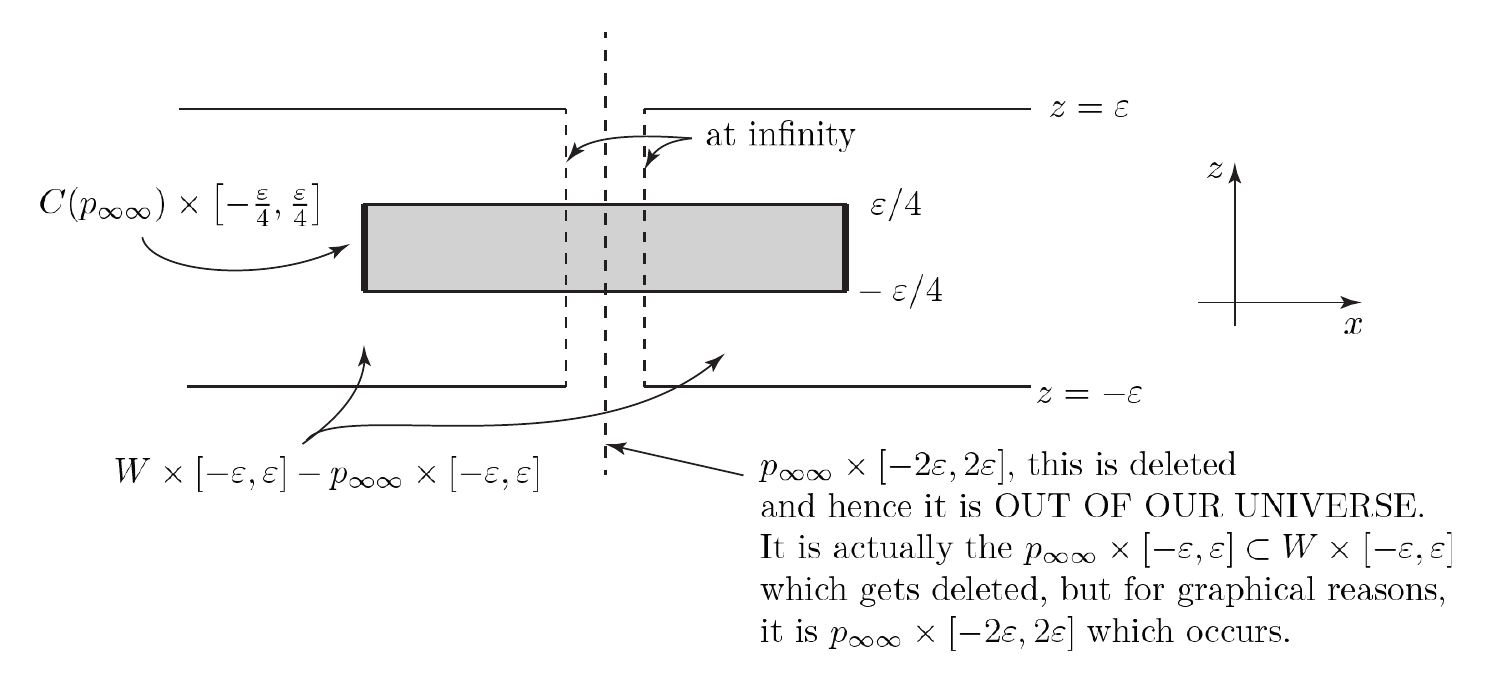}
$$
\label{fig6}
\centerline {\bf Figure 6.} 
\begin{quote}
We see here a section through $\Theta^3 (fX^2)$. The $W$ is the square in Figure 3. The compensating handle of index $\lambda = 2$ is the shaded $D^2 (p_{\infty\infty}) \times \left[ - \frac{\varepsilon}4 , \frac{\varepsilon}4 \right]$. The $C(p_{\infty\infty}) \times \left[ - \frac{\varepsilon}4 , \frac{\varepsilon}4 \right]$ $\left(\underset{\rm TOP}{=} S^1 \times [0,1]\right)$ is a singularity of $\Theta^3 (fX^2)$, {\bf not} of the undrawable type, but quite related. We will never smoothen it.
\end{quote}

\smallskip

Notwithstanding the fact that we try to keep the present survey of reasonable length and at reasonable distance from heavier technicalities, we have to give now 

\medskip

\noindent {\bf Some more details concerning $\Theta^3 (fX^2)$.}

\smallskip

To begin with, something like the Figure 6 might refer both to the $p_{\infty\infty}$ (proper) from Figure 3 and to the more difficult $p_{\infty\infty} (S)$'s. In both cases the vertical line $p_{\infty\infty} \times [-\varepsilon , \varepsilon]$ is deleted and compensated by the 2-handle $D^2 (p_{\infty\infty}) \times \left[-\frac{\varepsilon}2 , \frac{\varepsilon}2 \right]$, with the attaching zone $C(p_{\infty\infty}) \times \left[-\frac{\varepsilon}4 , \frac{\varepsilon}4 \right]$. But, while for $p_{\infty\infty}$ (proper) our Figure 6, as it stands, is quite realistic, it is far from being so, in the case of the $p_{\infty\infty} (S)$'s. In that case, to our Figure 6 one has to add the following feature, which I felt unable to draw. Each $p_{\infty\infty} (S) \in X^2$ is an accumulation point of immortal singularities of $X^2$,  all of the undrawable type.

\smallskip

This having been said, there are two versions for the $3^{\rm d}$ thickening of $fX^2$, the singular $\Theta^3 (fX^2)$ and the even more singular $\Theta^3 (fX^2)'$. They are both locally-finite.

\smallskip

Strictly speaking, what we have been discussing, so far, was $\Theta^3 (fX^2)'$, where ALL the $p_{\infty\infty} \times [-\varepsilon , \varepsilon]$ are deleted and compensated by 2-handles $D^2 (p_{\infty\infty}) \times \left[-\frac{\varepsilon}4 , \frac{\varepsilon}4 \right]$.

\smallskip

Now, in the case $p_{\infty\infty}$ (proper) ({\bf but} not in the case $p_{\infty\infty} (S)$) the Hole $p_{\infty\infty} \times [-\varepsilon , \varepsilon]$ can be smoothly HEALED by a process which I will explain (quite schematically) now.

\smallskip

Start from Figure 6 as it stands, with $p_{\infty\infty} \times [-\varepsilon , \varepsilon]$ and the handle $D^2 (p_{\infty\infty}) \times \left[-\frac{\varepsilon}4 , \frac{\varepsilon}4 \right]$ attached singularity. At the level of that figure, perform now the following steps:
\begin{enumerate}
\item[i)] Via an (infinite) $3^{\rm d}$ collapse which starts at $p_{\infty\infty} \times [-\varepsilon , \varepsilon]$ (living at infinity), delete enough material so that the attachement of the $2^{\rm d}$ handle should become smooth and not singular any longer.
\item[ii)] Then, via an (infinite) $3^{\rm d}$ dilatation, change the result of step i) into $W \times [-\varepsilon , \varepsilon] - p_{\infty\infty} \times \{\pm\varepsilon\}$, with the 2-handle smoothly melted into it. All this is now smooth.
\end{enumerate}

\smallskip

Globally, when one does this for all the $p_{\infty\infty}$ (proper), this defines a transformation
$$
\Theta^3 (fX^2)' \overset{{\rm HEALING}}{=\!\!=\!\!=\!\!=\!\!=\!\!=\!\!=\!\!=\!\!=\!\!=\!\!=\!\!\Longrightarrow} \{ \Theta^3 (fX^2), 
$$
$$
\mbox{an object where now only the $p_{\infty\infty} (S) \times [-\varepsilon , \varepsilon])$ stay deleted $+$ compensated} \}.
$$
For these twin $3^{\rm d}$ objects we will use the notation $\Theta^3 (fX^2)^{(')}$, meaning ``$\Theta^3 (fX^2)'$ OR $\Theta^3 (fX^2)$''.

\smallskip

And, with this, we open another subsection, namely:


\subsection{A toy model}

\medskip

The situation we will present now is highly simplified from several respects: there are no immortal singula\-rities (\`a la ${\rm Sing} \, \widetilde M^3(\Gamma)$), there are no triple points and no $p_{\infty\infty}$'s either.

\smallskip

But then there will not be any group action present now either. We replace (1.6) by the following map which {\ibf is} our toy model
$$
X^2 \overset{f}{-\!\!\!-\!\!\!-\!\!\!-\!\!\!\longrightarrow} R^3 \, , \leqno (1.16)
$$
where

\medskip

\noindent (1.16.1) \quad $X^2 \equiv \{$the $2^{\rm d}$ region called $R$, homeomorphic to $I \times {\rm int} \, I$ and PROPERLY embedded inside $(z=0) \subset R^3 = \{x,y,z \}$, which is displayed in the figure~7-(A)$\} \cup \underset{n=1}{\overset{\infty}{\sum}} D_n^2$, with $R$ and each of the $D_n^2$'s glued along $[\sigma_n , \Sigma_n] + [s_n , S_n]$. Here $D_n^2 = \{$the disk of diameter $[\sigma_n , s_n] = R \cap (x=x_n)$, living in the plane $(x=x_n) \subset R^3$ and displayed in the figure~7-(B)$\}$.

\medskip

Figure~7 should explain how the spare parts of $X^2$, listed in (1.16.1) are glued together. Restricted to each of these spare parts, the map $f$ is the natural embedding into $R^3$.

$$
\includegraphics[width=150mm]{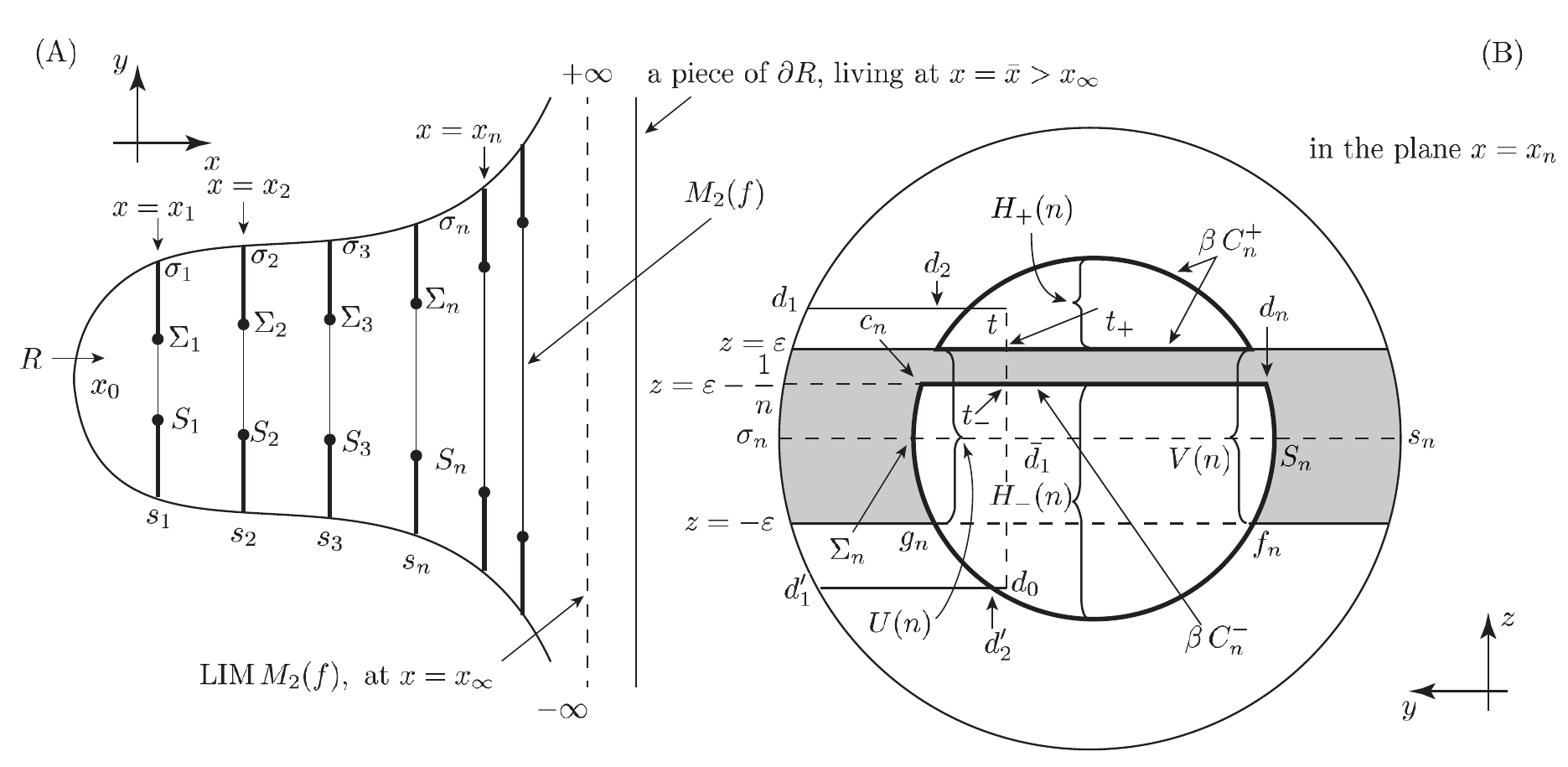}
$$
\label{fig7}
\centerline {\bf Figure 7.} 

\smallskip

\begin{quote} 
We display here the spare parts $R$ and $D_n^2$ of $X^2$. The $\Sigma_n$, $S_n$ are mortal singularities of the $f$ (1.16), and the curvilinear arcs $U(n) \ni \Sigma_n$, $V(n) \ni S_n$ stand for undrawable immortal singularities of $\Theta^3 (X^2)$. The $H_{\pm} (n)$ are Holes, open sets, drilled out of $D_n^2$.
\end{quote}

\smallskip

We will define the following smooth non-compact manifold with large boundary, to be compared to the (1.15) above,
$$
\Theta^3 (fX^2) \equiv \left( R \times [- \varepsilon \leq z \leq \varepsilon] \right) \cup \sum_1^{\infty} D_n^2 \times [- \varepsilon_n \leq x \leq \varepsilon_n] - {\rm LIM} \, M_2(f) \times \{z = \pm \, \varepsilon\} \, . \leqno (1.17)
$$
We also introduce $\Theta^{N+4} (X^2)$ and $\Theta^{N+4} (fX^2) \equiv \{\Theta^4$ (of the already smooth $\Theta^3 (fX^2)) = \Theta^3 (fX^2) \times I\} \times B^N$. In preparation for what we will do afterwards in real life (no longer in a toy model context), we will consider now a transformation $\Theta^{N+4} (X^2) \Rightarrow \Theta^{N+4} (fX^2)$, a high-dimensional geometric version of the zipping process of $f$, which should consist of an infinite sequence of steps, each of which, individually, should in principle at least be GSC preserving.

\smallskip

Notice, for instance, that the deletion of the ${\rm LIM} \, M_2(f)$ contribution in formula (1.17) can also be described as an infinite dilatation, sending the bands ${\rm LIM} \, M_2(f) \times \{ z = \pm \, \varepsilon \} \times [x_{\infty} , x_{\infty} + \varepsilon_0]$ for a very small $\varepsilon_0$ to infinity, a GSC preserving step.

\smallskip

In real life, our object of interest will be a cell-complex $S_u (\widetilde M^3 (G))$, a distant analogue of the present $\Theta^{N+4} (fX^2)$. But, essentially because of the $p_{\infty\infty}$'s, $S_u \, (\widetilde M^3 (G))$ is just a cell-complex and not a smooth manifold. On the other hand, the main step in the proof that $G \in {\rm QSF}$ will be to show that $S_u (\widetilde M^3 (G)) \in {\rm GSC}$. In that proof, the group action of $G$ (``{\bf discrete symmetry with compact fundamental domain}'') will play an essential role.

$$
\includegraphics[width=160mm]{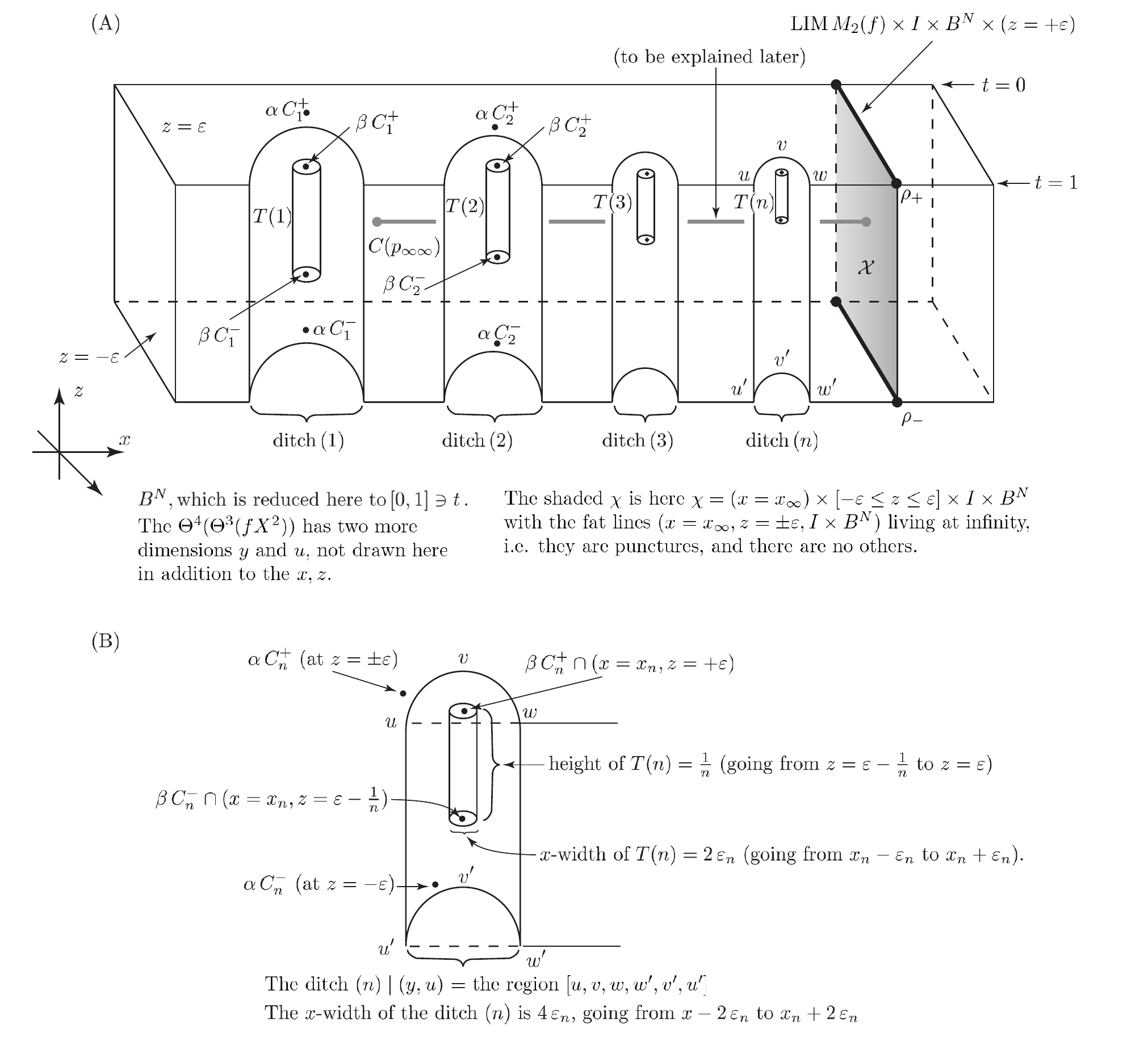}
$$
\label{fig8}
\centerline {\bf Figure 8.} 

\smallskip

\begin{quote} 
Everything we see both in (A) and in (B), lives at some generic values $y \in [y \, V(n),y \, U(n)]$ and $u \in I$ (fourth coordinate in $\Theta^4$), with $V(n) , U(n)$ like in figure~7-(B). In the present (B) we see the ditch $(n) \mid y$ and, with this, ditch $(n) \equiv \bigcup \, {\rm ditch} \ (n) \mid (y,u)$, with $(y ([t_- , t_+]$ of $V(n)) \leq y \leq (y ([t_- , t_+]$ of $U(n))$. Here $[t_- , t_+]$ is like in figure 7-(B), i.e. $(y([t_- , t_+]$ of $V(n))$ (respectively $(y([t_- , t_+]$ of $U(n)))$ is slightly larger than $yV(n)$ (respectively slightly smaller than $y \, U (n)$); see the $y$ coordinates in the Figure 7-(B).
\end{quote}

\bigskip

In our toy-model, there is no group $G$ and the $\Theta^{N+4} (fX^2)$ will possibly fail to be GSC. But we will do, starting from the $\Theta^{N+4} (X^2) \in {\rm GSC}$ the correct moves, namely TURNING QUOTIENT SPACE PROJECTIONS INTO INCLUSION MAPS, (which is the key idea of what we call the high-dimensional geometric realization of the zipping process). These moves do go in the good general GSC direction. It will be pedagogically useful to see how that goes in the context of the toy-model transformation $\Theta^{N+4} (X^2) \Rightarrow \Theta^{N+4} (fX^2)$ and nevertheless we fail to get the kind of PROPER homotopy of links, (which in real life will be a key to $S_u (\widetilde M^3 (G))$ being QSF), in our present case. So, just watch what is coming next, and keep it then in mind when we will go to real life.

\smallskip

For our transformation we will introduce three ingredients, to be very much used, afterwards, in real life too. These are:
$$
\mbox{Holes, Ditches and {\bf Partial} ditch-filling.} \leqno (1.18)
$$
The Holes do get eventually compensated by 2-handles, of course. In the Figure 7-(B) we have presented the Holes $H_{\pm} (n) \subset D_n^2$. To be precise, in Figure 7-(B) we see a smaller disc $d_n^2 \subset D_n^2$, concentric to $D_n^2$, with
$$
\partial \, d_n^2 \supset \{ U (n) , V(n)\} \, ,
$$
and with $H_{\pm} (n) \subset d_n^2$. Most of the fat curves $\beta \, C_n^{\pm}$, lives in $\partial d_n^2$. In the Figure 7-(B), for the $H_{\pm} (n)$'s the boundaries are in fat lines and their interiors get deleted. (Unlike what happens with the deletion of the $p_{\infty\infty}$'s, when closed sets are getting deleted, open sets are deleted now.) In our figure, the shaded areas which are {\bf outside} of the central $d_n^2$ (inside which $H_{\pm} (n)$ live), correspond to identifications at level of $\Theta^3 (X^2)$, leaving us with the singularities $U(n)$, $V(n)$.

\smallskip

With $C_n^{\pm} \equiv \partial \, {\rm Hole} \, H_{\pm} (n)$, comes a framed link
$$
\left\{\mbox{framed link} \ \sum_n C_n^{\pm} \right\} \xrightarrow{ \ \beta \ } \partial \, \Theta^{N+4} (X^2 - H) \, , \leqno (1.19)
$$
where $\beta$ is like in Figure 7-(B), and ``$-H$'' means ``with the Holes deleted''.

\smallskip

The shaded area contained in the inside of the central $d_n^2$, between $\varepsilon - \frac1n \leq z \leq \varepsilon$, corresponds to a piece of  $(D_n^2 - H) \cap f^{-1} \left(R \times \left[\varepsilon \geq z \geq \varepsilon - \frac1n \right] \right)$. It stays far from the identifications at level $\Theta^{N+4} (X^2 - H)$, but it is part of the identifications leading to $\Theta^{N+4} (fX^2-H)$. And now we have a second framed link
$$
\left\{\mbox{framed link} \ \sum_n C_n^{\pm} \right\} \xrightarrow{ \ \alpha \ } \partial \, \Theta^{N+4} (fX^2 - H) \, . \leqno (1.20)
$$
Here $\alpha \, C_n^+ = \beta \, C_n^+$, but unlike the $\beta \, C_n^-$ which makes it all the way up to $z = \varepsilon - \frac1n$ the $\alpha \, C_n^-$ only makes it to $z = -\varepsilon$. At the level $\Theta^{N+4} (fX^2 - H)$, the corresponding $H_- (n)$ stops at $z = - \varepsilon$. The two figures 8 may explain the positions of $\alpha \, C_n^{\pm}$, $\beta \, C_n^{\pm}$. This ends the story of the Holes in the toy model.

\smallskip

The Ditches are suggested in Figure 8. The ditches are indentations made inside $\Theta^{N+4} (X^2 - H)$, {\bf making use of the supplementary dimensions $B^N$}, which in Figure 8 are reduced to a $B^N = [0,1]$, for typographical reasons. And what we suggest in Figure 8 is actually
$$
\sum_{n=1}^{\infty} {\rm ditch} \, (n) \, \vert (y,u) \subset \Theta^{N+4} (X^2-H) \vert \, R \, .
$$
Here $y$ is a fixed value of the third dimension of $\Theta^3 (X^2)$, and in terms of the Figure 7, the spectrum of values of $y$ is exactly the following
$$
y (V(n)) \geq y \geq y (U(n)) \, , \leqno (1.21)
$$
while $u$ is any generic point of $I \ni u$, fourth coordinate of $\Theta^4 = \Theta^3 \times I$, fixed too. Our ${\rm ditch} \, (n) \mid (y,u)$ is concentrated around the arc $A (x_n , y,u) = (x = x_n , y , -\varepsilon \leq z \leq \varepsilon , u , \{$some point $p \in \partial \, B^N$, and in Figure 8 this is $t=1\}) \, \subset$
$$
\subset \ (\underbrace{R \, (\mbox{Figure 7-(A)}) \times [-\varepsilon , \varepsilon] \times I}_{{\rm in} \ \Theta^4} \times B^N) \mid (y,u) \, .
$$

Then ${\rm ditch} \, (n) \mid (y , u) = A(x_n , y , u) \times b^N (n)$, with $b^N (n)$ an $N$-ball $b^N(n) \subset B^N$, concentrated around $p \in \partial \, B^N$ ($t=1$, in the context of Figure 8), and coming with
$$
\lim_{n \, = \, \infty} {\rm diam} \, b^N (n) = 0 \, ,
$$
and, finally
$$
{\rm ditch} \, (n) = \bigcup_{\overbrace{\mbox{\footnotesize$y \in (1.21) , \, u \in I$}}} {\rm ditch} \, (n) \mid (y,u) \, . \leqno (1.22)
$$

Figure 8-(B), which focuses on a detail of 8-(A), presents us also with a piece of $D_n^2-H$, namely the
$$
T(n) \mid (y,u) = \left( \Theta^{N+4} (X^2-H) \left\vert \left\{ (D_n^2 -H) \cap \left[y \in (1.21) , z \in \left[ \varepsilon - \frac1n , \varepsilon \right] \right] \right\} \right)\right\vert (y,u) \, . \leqno (1.23)
$$

The $\underset{\overbrace{\mbox{\footnotesize$y \in (1.21) , \, u \in I$}}}{\bigcup} T(n) \mid (y,u)$ is a $4^{\rm d}$ 1-handle $T(n)$ and from Figure 8 we can read a SPLITTING
$$
\Theta^{N+4} (X^2-H) \ \mbox{(which is non-singular)} - \{{\rm DITCHES}\}  \, =
$$
\vglue -8mm
$$
= \  \{{\rm MAIN} \ \Theta^{N+4} (X^2-H)\} + \sum_{n=1}^{\infty} \{ (N+4)-{\rm handle} \ T(n) \ \mbox{of index} \ \lambda = 1 \} \, . \leqno (1.24)
$$
Here, essentially, we have $\{{\rm MAIN} \ \Theta^{N+4} (X^2-H)\} \approx R \cup \underset{n \, = \, 1}{\overset{\infty}{\sum}} \, (D_n^2 - d_n^2)$.

\smallskip

There is a PROPER embedding $j$ entering in the following commutative diagram
$$
\xymatrix{
\Theta^{N+4} (X^2-H) - \{{\rm DITCHES}\} \ar@{->>}[d] \ar[rr]_-j &&\Theta^{N+4} (fX^2 - H) \ar@{->>}[d] \, . \\
X^2 - H \ar[rr]_-f &&fX^2-H
}
\leqno (1.25)
$$
Here $j \mid (T(n) \mid (y,u))$ sends $T(n) \mid (y,u)$ inside the following ditch
$$
{\rm ditch} \, (n) \mid (y,u) - \partial ({\rm ditch} \, (n) \mid (y,u)) \subset \{{\rm DITCHES}\} \mid (y,u) \subset \Theta^{N+4} (fX^2-H) \mid (y,u) \leqno (1.26)
$$
in the manner which is suggested in the Figure 8-(B).

\smallskip

This is the story of the ditches in the context of the Toy-model. As the Figure 8-(B) suggests, the $j$-image of $T(n)$ really goes in the void left behind by the deleted DITCH, now just an indentation at level $B^N$ (supplementary dimension). So we do have $j \, T(n) \subset {\rm ditch} \, (n) \subset [-\varepsilon \leq z \leq + \varepsilon]$ but, actually, both the $j \, T(n)$, as well as rest of the ditch-filling material, a glue which will later also be send in the ditch, so as to join the 
$j \, T(n)$ to the lateral surface of the ditches will always be confined to the $\left[ \varepsilon - \frac1n , \varepsilon \right] \subset [-\varepsilon , \varepsilon]$.

\smallskip

Here
$$
\lim_{n \, = \, \infty} \left[\varepsilon - \frac1n \leq z \leq \varepsilon \right] \cap A(x_n , y , u)
$$
$$
\in \ \{{\rm LIM} \, M_2 (f) \times I \times B^N \times (z = + \varepsilon) , \ \mbox{which is {\bf deleted} at the level of} \ \Theta^{N+4} (fX^2-H)\} \, ,
$$
and our $\underset{n \, = \, \infty}{\lim}$ is just a point, let us say the $\rho_+$ in Figure 8-(A). Another way to phrase things, is that we have
$$
\lim_{n \, = \, \infty} T(n) \subset {\rm LIM} \, M_2 (f) \times I \times B^N \times \{ z = + \varepsilon \} \, ,
$$
and the RHS of this formula lives at infinity.

\smallskip

This is the phenomenon of {\bf partial ditch filling}, which goes hand in hand with the position of the Hole $H_- (n)$ in Figure 7-(B). Finally, the next lemma expresses what is achieved at the Toy-model level, and what is not.

\bigskip

\noindent {\bf Lemma 6.} 1) {\it Inside $\Theta^{N+4} (fX^2-H)$ as ambient space, we can perform an infinite system of smooth dilatations and additions of handles of index $\lambda > 1$, call this infinite process

\bigskip

\noindent {\rm (1.27)} \quad $j (\Theta^{N+4} (X^2-H) - \{\mbox{DITCHES}\}) \underset{{{\rm DIL} \atop \mbox{\footnotesize (like ``dilatation'')}}}{=\!\!=\!\!=\!\!=\!\!=\!\!=\!\!=\!\!=\!\!=\!\!=\!\!=\!\!=\!\!=\!\!=\!\!\Rightarrow} \{$the $j (\Theta^{N+4} (X^2-H) - \{\mbox{DITCHES}\})$ with all that part of the $\underset{n \, = \, 1}{\overset{\infty}{\sum}}$ ditch $(n)$ which is not yet occupied by $j \ \underset{n}{\sum} \ T(n)$ FILLED IN, but ONLY PARTIALLY, at level $\varepsilon - \frac1n \leq z \leq \varepsilon \}$.}

\medskip

2) {\it The step DIL above is PROPER inside the ambient space $\Theta^{N+4} (fX^2-H)$. We will denote by
$$
S_b (fX^2-H) \equiv \{\mbox{the target $\{\ldots\}$ of DIL, in the RHS of the formula {\rm (1.27)} above}\} \, . \leqno (1.28)
$$

This $S_b (fX^2-H)$, living at the target of our transformation DIL, is a smooth $(N+4)$-dimensional manifold. The reason for our notation ``$S_b$'' will be explained later, in due time. For further purposes, one should notice that all the individual steps which occur inside DIL are of the GSC-preserving type.

\smallskip

The $S_b$ will occur big, again, when we will move from toy-model to real life.}

\medskip

3) {\it There is a commutative diagram of PROPER embeddings
$$
\xymatrix{
\Theta^{N+4} (X^2-H)-\{\mbox{DITCHES}\} \ar[rr]^-j \ar[d]^{\mbox{DIL}} &&\Theta^{N+4} (fX^2-H) \, . \\
S_b (fX^2-H) \ar[urr]_{\mathcal J}
}
\leqno (1.29)
$$
}

4) {\it There is, also a {\ibf simple-minded diffeomorphism} $\eta$, which enters into the diagram below
$$
\xymatrix{
&&\Theta^{N+3} (\partial X^2) \ar[dl]^b \ar[dr]^a \\
&S_b (fX^2-H) \ar[rr]^{\eta} &&\Theta^{N+4} (fX^2-H) \, , \\
&&\underset{n}{\sum} \, C_n^{\pm} \ar[ul]^{\beta}\ar[ur]_{\alpha}
}
\leqno (1.30)
$$
where $\partial X^2 = \partial R \cup \underset{n \, = \, 1}{\overset{\infty}{\sum}} \partial D_n^2$ (see Figure $7$) and $\Theta^{N+3} (\partial X^2)$ is its obvious $N+3$ smooth thickening, where $a,b$ are the obvious inclusions and $\alpha , \beta$ are like in {\rm (1.20), (1.19)} respectively. Also
$$
{\rm Im} \, b \cap {\rm Im} \, \beta = \emptyset = {\rm Im} \, a \cap {\rm Im} \, \alpha \, .
$$

Here is what we can say concerning the diagram {\rm (1.30)}:
\begin{enumerate}
\item[$\bullet$)] The upper triangle and the restriction of the lower triangle to $\underset{n}{\sum} \, C_n^+$ commute, strictly.
\item[$\bullet$$\bullet$)] The embedding ${\mathcal J}$ in {\rm (1.29)} is connected to $\eta$ by a {\ibf not} boundary respecting isotopy.
\item[$\bullet$$\bullet$$\bullet$)] The lower triangle, when restricted to $\underset{n}{\sum} \, C_n^-$ also commutes, but only {\ibf up to homotopy}.
\end{enumerate}
}

\medskip

5) {\it We move now to the smooth $(N+4)$-manifolds

\bigskip

\noindent {\rm (1.30.1)} \quad $\Theta^{N+4} (fX^2) = \Theta^{N+4} (fX^2-H) + \underset{n}{\sum} \, \{$the $(N+4)$-dimensional $2$-handles $D^2 (\alpha \, C_n^{\pm})\}$, which in view of what will happen, later on, in real life, we also call $S_u (fX^2)$, and to $S_b (fX^2) \equiv S_b (fX^2-H) + \underset{n}{\sum} \, \{$the $(N+4)$-dimensional $2$-handles $D^2 (\beta \, C_n^{\pm})\}$.

\bigskip

With this, there is NO homeomorphism of the form
$$
\mbox{\hglue -15mm} (S_b (fX^2) , \Theta^{N+3} (\partial X^2)) \longrightarrow (\Theta^{N+4} (fX^2) , \Theta^{N+3} (\partial X^2)).
$$
\vglue-7mm
$$
{\mid\mbox{\hglue -2mm}}_{\overset{\rm id}{-\!\!-\!\!-\!\!-\!\!-\!\!-\!\!-\!\!-\!\!-\!\!-\!\!-\!\!-\!\!-\!\!-\!\!-\!\!-\!\!-\!\!-\!\!-\!\!-\!\!-\!\!-\!\!-\!\!-\!\!-\!\!-\!\!-\!\!-\!\!-\!\!-\!\!-\!\!-\!\!-\!\!-\!\!-\!\!-\!\!-\!\!-\!\!-\!\!-\!\!-\!\!-\!\!-\!\!-\!\!-\!\!-}}{\mbox{\hglue -2mm}\uparrow}
$$
It follows that no such diffeomorphism exists either, and hence the lower triangle in {\rm (1.30)} does NOT commute up to PROPER homotopy.}

\bigskip

The last item in Lemma 6 shows that the toy-model fails to deliver an ingredient which, in real life, will be essential, namely commutativity up to PROPER homotopy. But there is a way to redeem this, which I will briefly describe now. Our $\Theta^3 (fX^2)$ from (1.17) (which, of course, should not be mixed up with the real-life $\Theta^3 (fX^2)$), contains the following BAD RECTANGLE
$$
{\mathcal R}_{\infty} = \{ x=x_{\infty} , - \infty < y < \infty , - \varepsilon \leq z \leq \varepsilon \} \subset R \times [-\varepsilon,\varepsilon] \, ,
$$
and we use here the notations from (1.16.1).

\smallskip

This bad rectangle is to be sent to infinity, {\bf not} by deletion, {\bf but} by $3^{\rm d}$ (infinite) dilatation, namely by adding ${\mathcal R}_{\infty} \times [0,\infty)$ along ${\mathcal R}_{\infty} = {\mathcal R}_{\infty} \times \{0\}$, at the level of $\Theta^{N+4} (fX^2 - H)$. When one does this, the lower triangle in (1.30) commutes now up to PROPER homotopy.

\smallskip

So, we end here the discussion of the toy model and turn back to real life, but keep in mind that the various items introduced with the toy model, like (1.18) and the whole discussion from Lemma 6 will reverberate through the real life story. From now on, $(X^2 , f)$ will be like in the Theorem 5 and no longer like in the toy-model.

\smallskip

We move now to the final and most important subsection of the present Section 1, namely

\subsection{The geometric realization of the zipping in high dimensions}

In the present subsection, which is the core of the matter of the present survey, we will introduce the four functors $S_u , S'_u , S_b , S'_b$, which are the real-life analogues of the $\Theta^{N+4}$, $S_b$ from the toy-model. (Here, ``$S$'' should mean ``structure'', while the subscripts ``$u$'', ``$b$'' stand respectively, for ``usual'' and ``bizarre''.)

\smallskip

We will then present the plan of the proof that $S_u (\widetilde M^3 (\Gamma))$ is GSC, which is the main ingredient for the proof of Theorem 3 (``all $G$'s are QSF'').

\smallskip

When we move to real life, there is our group $G$ and everything is equivariant too. We will find there again, the analogue of the lower triangle in (1.30), with objects $S_b$ and $S_u$ to be properly introduced later, which are connected now to $\widetilde M^3 (G)$ and not to $\Theta^3 (fX^2)$ from the toy-model.

\smallskip

Then we can quotient by $G$ and move downstairs at level $M^3 (G)$, making use of the functoriality of the $S_b$, $S_u$ (appearing now in lieu of the $\Theta^{N+4} (fX^2-H)$). But even then, the $S_b (M^3 (G)-H)$, $S_u (M^3(G)-H)$ will fail to be compact and the $\{ C_n^{\pm}\}$ will continue to be infinite.

\smallskip

Nevertheless, the fact that $M^3 (G)$ is compact, combined with the uniform boundedness of the zipping flow (see Theorem 4 in Part I), will make that, at level $M^3 (G)$, the analogue of the lower triangle in (1.30) commutes up to PROPER homotopy. And, in the range of dimensions which concern us, PROPER homotopy implies isotopy, and that is what we need, in real life.

\smallskip

In the beginning of the present section we have mentioned that the main ingredient of the proof of Theorem 3, the object of the Parts [II], [III] of the trilogy is the proof that a certain $(N+4)$-dimensional cell-complex $S_u (\widetilde M^3 (G))$ is GSC. Its exact definition is quite complicated but by now we can say this much: Our real-life $S_u (\widetilde M^3 (G))$ plays the role of $\{\Theta^{N+4} (fX^2)$ from (1.30.1)$\}$, of which it is the real life analogue.

\smallskip

Of course, our real life (1.6) is much more complicated than the $(X^2,f)$ of the toy model, so a lot of embellishments with respect to (1.30.1) are necessary, making in particular that $S_u (\widetilde M^3 (\Gamma))$ cannot be a smooth manifold, but just a cell-complex with very controlled singularities.

\smallskip

But right now we would like to explain the reason behind the {\bf partial} ditch-filling, something which will be essential for the proof of $S_u (\widetilde M^3 (\Gamma)) \in {\rm GSC}$.

\smallskip

The process DIL, introduced in (1.27), is an attempt of geometric high-dimensional realization of the zipping process (in the toy-model context) i.e. the replacement of the big quotient-space projection $X^2-H \to fX^2-H$, via a sequence of inclusion maps (which are much more manageable than quotient-space projection), followed by $2$-handle additions like in (1.30.1). The inclusion maps mentioned here are, of course, the process of PARTIAL DITCH FILLING.

\smallskip

The next Lemma 7 is certainly NOT part of the proof of Theorem 3, but it is an illustration which we present here for pedagogical purposes, of what could go wrong in the context of our present discussion, if we would proceed with a complete ditch-filling and not with a partial one, as we actually do, with its $\underset{n \, = \, \infty}{\lim} \, T(n) = \infty$.

\bigskip

\noindent {\bf Lemma 7.} (Po-Tanasi \cite{41}.) 1) {\it Let $V^3$ be an open simply-connected $3$-manifold, without any other additional conditions, and let also $m \in Z_+$ be high enough. There is then a sequence of smooth non-compact $(m+3)$-manifolds, with very large boundaries, and of smooth inclusions
$$
X_1^{m+3} \subset X_2^{m+3} \subset X_3^{m+3} \subset \ldots
$$
such that: {\rm 1)} $X_1^{m+3} \in \mbox{GSC}$ and all the inclusions above are finite combinations of smooth Whitehead dilatations and/or additions of handles of index $\lambda > 1$, hence GSC-preserving moves.}

\medskip

2) {\it Let $\varinjlim X_i^{m+3}$ be the union of the objects above, endowed with the {\ibf weak topology} (and there is no other reasonable one which one can use in our context); then there is a {\ibf continuous bijection}
$$
\varinjlim X_i^{m+3} \xrightarrow{ \ \psi \ } V^3 \times B^m \, .
$$
}

Remember here that in the weak topology, which is considered, a set $F \subset \varinjlim X_i^{m+3}$ is closed iff all the $F \cap X_i^{m+3}$'s are closed. But then, I claim that the $\psi^{-1}$ cannot be continuous. Here is the proof. If $\psi$ would be a homeomorphism then $V^3 \times B^m$ would be a GSC manifold. But then, also, our $V^3$ might well be the Whitehead manifold ${\rm Wh}^3$ which comes with $\pi_1^{\infty} \, {\rm Wh}^3 \ne 0$. This would contradict the stabilization lemma from part I of this survey. (And do not mix up here the ``$I$'' in this survey and the ``$I$'' in the trilogy.)

\smallskip

With these things, it follows from an old classical theorem of Brower, that, $\varinjlim X_i^{m+3}$ cannot be a manifold. Even more, and this is the point here, our $\varinjlim X_i^{m+3}$ is not a metrizable space. We will call this the non-metrizability barrier, a typical difficulty encountered when trying to prove Theorem 3. The {\bf partial} ditch-filling is necessary for overcoming this barrier. There are also other barriers, like the so-called {\bf Stallings barrier}, to be discussed later when we will come to discuss Part III of the Trilogy, in the next section.

\smallskip

Then there is also the following barrier. On our $X^2$ in (1.6) there are two not everywhere well-defined flows, the collapsing flow, connected to $X^2$ being GSC, and the zipping flow. Individually, each of these two flows is quite nice and tame, but we have transversal contacts for the two kind of trajectories
$$
\{\mbox{zipping flow}\} \pitchfork \{\mbox{collapsing flow}\} \, ,
$$
and this can create close oriented loops, which are highly dangerous for us.

\smallskip

We will discuss now, a bit, the real-life $S^{(')}_u (\widetilde M^3 (G))$. Corresponding to the two versions for the $3^{\rm d}$ thickening $\Theta^3 (fX^2)$ and $\Theta^3 (fX^2)'$, there are two versions for the $(N+4)$-dimensional complex $S_u (\widetilde M^3 (G))$ and $S'_u (\widetilde M^3 (G))$. The $S'_u$ is more singular than $S_u$, and $S^{(')}_u$ will mean ``$S_u$ or $S'_u$''. Constructing the $S_u^{(')} (\widetilde M^3 (G))$ is a process which happens by successive stages, where more ingredients and more dimensions are, one by one, thrown in and we will give here a schematical view of its structure. The first step is to consider, in the world of Theorem 5, the $\Theta^3 (fX^2)^{(')}$, and we had already a glimpse of this locally finite singular, two version  analogue of the smooth $\{ \Theta^3 (fX^2)$ from the context of our Toy-model, see here the Figures 7, 8$\}$.

\smallskip

Then, as in (1.10), we introduce the

\bigskip

\noindent (1.31) \quad $\Sigma_1 (\infty) \subset \widetilde M^3 (G)$, the union of the limiting positions of the compact walls $W$(BLUE), $W$(RED), $W$(BLACK). But these last ones are mute in (1.11.1), (1.11.2).

\bigskip

We have the $2^{\rm d}$ non-compact object $\Sigma (\infty) \equiv \{\Theta^3 (fX^2) \cap \Sigma_1 (\infty)$ with the boundary $(\partial \Sigma (\infty))$ deleted$\}$, a formula which should be compared to (1.11.2) of which it is the $3^{\rm d}$ version. The $\Sigma (\infty)$ is an open $\{$surface with branching lines, like $\{$figure $Y\} \times R\}$, and it is a $(-\varepsilon , + \varepsilon)$-thickening of $f \, {\rm LIM} \, M_2 (f)$ (1.11.2), with a boundary living at $\{\pm \, \varepsilon\}$ i.e. at infinity, and not at finite distance in the real world. It will be denoted by $\partial \, \Sigma (\infty)$. It corresponds to the deletions performed in (1.15) and (1.17). In Figure 5, $\partial \, \Sigma (\infty)$ occurs as the fat red points with the dotted red line in between standing for the open $\Sigma (\infty)$. We also have the branching locus of $\Sigma (\infty)$, namely
$$
\sum_{p(\infty\infty)} p(\infty\infty) \times (-\varepsilon , \varepsilon) + \sum_{p(\infty\infty)(S)} p_{\infty\infty} (S) \times (-\varepsilon , \varepsilon) \underset{\mbox{\footnotesize PROPER} \atop \mbox{\footnotesize embedding}}{\subset} \Sigma (\infty) \, , \leqno (1.31.A)
$$
which gets deleted from $\Theta^3 (fX^2)$ (and this is compensated by $2$-handles, like in Figure 6), before we can go further into the construction of $S_u^{(')} (\widetilde M^3 (G))$. And the next step is to send to infinity $\Sigma (\infty) - \underset{p_{\infty\infty {\rm (ALL)}}}{\Sigma} \, p_{\infty\infty} \times (-\varepsilon, \varepsilon) \equiv \Sigma_0 (\infty)$, via the non-compact dilatation below, and {\bf not} by any deletion, meaning the following thing

\bigskip

\noindent (1.31.B) \quad We add to $\Theta^3 (fX^2) - \underset{p_{\infty\infty}}{\Sigma} \, p_{\infty\infty} \times (-\varepsilon, \varepsilon)$ the $\Sigma_0 (\infty) \times [0,\infty)$, along $\Sigma_0 (\infty) \times \{0\} = \Sigma_0 (\infty)$.

We have already introduced the $\Theta^3 (fX^2)$ and $\Theta^3 (fX^2)'$ coming with $p_{\infty\infty} \times (-\varepsilon , \varepsilon)$'s deleted and compensated by $2$-handles $D^2 (p_{\infty\infty})$, and which are both locally-finite. With this, a first approximation version (forgetting various embellishments and subtelties) of $S_u^{(')} (\widetilde M^3 (G))$ is
$$
S_u^{(')} (\widetilde M^3 (G)) \equiv \Theta^4 (\Theta^3 (fX^2)',{\mathcal R} ) \times B^N \, .
$$

But, the real definition of $S_u^{(')}$ has to integrate the addition of $\Sigma_0 (\infty) \times [0,\infty)$ from (1.31.B), but how to do this correctly at the level of $(N+4)$-dimensions will not be explained here.

\smallskip

Retain, also, that all this constructions are equivariant, like in (1.31.A) below.

\smallskip

With this, the key step for the proof of Theorem 3, is the following result, on which we will concentrate in this section.

\bigskip

\noindent {\bf Theorem 8.} 1) {\it $S'_u (\widetilde M^3 (G))$ is GSC.}

\medskip

2) {\it There is a relatively easy transformation, proceeding on the lines of our healing process $\Theta^3 (fX^2)' \Rightarrow \Theta^3 (fX^2)$ mentioned earlier,
$$
S'_u (\widetilde M^3 (G)) \Longrightarrow S_u (\widetilde M^3 (G))
$$
and, with this, point {\rm 1)} implies that $S_u (\widetilde M^3 (G))$ is GSC too.}

\bigskip

A very essential ingredient for the proof of Theorem 8 is the fact that there exist $3^{\rm d}$ REPRESENTATIONS $X^3 \xrightarrow{ \ f \ } \widetilde M^3 (G)$ which are locally finite, equivariant, and with uniformly bounded zipping length. We will call them {\bf appropriate} REPRESENTATIONS. Then, we will want to be able to speak of $S_u , S'_u$ in the same breath and so we will write $S_u^{(')}$, when we want to encompass them both. And the big fact is that $S_u^{(')}$ is a FUNCTOR with good properties of localization and glueing. Without trying to make explicit our category-theory present context, the items below are useful illustrations of what our functoriality means.

\bigskip

\noindent (1.32.A) \quad We have free group actions
$$
G \times S_u^{(')} (\widetilde M^3 (G)) \to S_u^{(')} (\widetilde M^3 (G)) \, .
$$

\noindent (1.32.B) \quad Downstairs, at level $M^3(G)$, the $S_u^{(')} (M^3 (G))$'s are also well-defined.

\bigskip

With this, we have
$$
S_u^{(')} (\widetilde M^3 (G)) / G = S_u^{(')} ( M^3 (G)) \, .
$$
But this formula is not just a tautology. Of course, we may use it for {\it defining} the $S_u^{(')} (M^3 (G))$. But then, the point is that the $S_u^{(')} ( M^3(G))$ is also {\bf directly defined}, downstairs, by localizing the functor, and this is important. So, we have two distinct definitions for $S_u^{(')} (M^3 (G))$ and, very importantly, they are equivalent.

\bigskip

\noindent (1.32.C) \quad Again, by functoriality, we have
$$
S_u^{(')} (\widetilde M^3 (G))  = (S_u^{(')} ( M^3 (G)))^{\sim} \, .
$$

\noindent (1.32.D) \quad Everything in our little theory is equivariant. Like in (1.18) we will have holes $H$ (to which we will come back), ditches and partial ditch-filling, all of these items being equivariant too. Once the functor $S_u^{(')}$ has good localization properties, one can make sense of the $(N+4)$-dimensional $S_u^{(')} (\widetilde M^3 (G)-H)$ (with all holes, compact and open, deleted and not compensated) and, analogous to the first formula in (1.30.1) there is a {\bf reconstruction formula}
$$
S_u^{(')} (\widetilde M^3 (G)) = S_u^{(')} (\widetilde M^3 (G)-H) + \sum_H D^2 (H) \, ,
$$
where $D^2(H)$ are $(N+4)$-dimensional handles of index two, compensating for the holes $H$. Their attachments may be singular. To be precise, we have a framed link, for $\partial D^2 (H) = C(H)$,
$$
\sum_{\mbox{\footnotesize all $H$'s}} C(H) \underset{\mbox{\footnotesize PROPER embedding}}{\overset{\alpha}{-\!\!\!-\!\!\!-\!\!\!-\!\!\!-\!\!\!-\!\!\!-\!\!\!-\!\!\!-\!\!\!-\!\!\!-\!\!\!-\!\!\!-\!\!\!-\!\!\!-\!\!\!\longrightarrow}} S_u^{(')} (\widetilde M^3 (G)-H) \, ,
$$
via which the 2-handles get attached in the reconstruction formula. The $\alpha \, \underset{H}{\sum} \, C(H)$ is not confined into something like $\partial \, S_u^{(')} (\widetilde M^3 (G)-H)$, we are singular here, remember.

\smallskip

Moreover, {\bf reconstruction formulae} like the one just written down, function also downstairs, at the level of the $M^3 (G)$. End of (1.32.D).

\bigskip

Now, if $S_u^{(')} (M^3 (G))$ would be a compact object, then Theorem 8, together with the items like (1.32.A), (1.32.B), $\ldots$ above, would immediately imply that $G \in {\rm QSF}$. But the $S_u^{(')} ( M^3 (G))$ is NOT compact, although the $M^3(G)$ is. It is actually an infinitely foamy, $(N+4)$-dimensional version of the compact $M^3 (G)$.

\bigskip

\noindent [Terminology. The ``$S_u$'' stands for ``usual structure'', while the forthcoming ``$S_b$'' (the real-life analogue of (1.28) and of the second formula in (1.30.1)) stands for ``bizarre structure''.]

\bigskip

Finally, here is a last virtue which we want $S_u^{(')}$ to have, of a different kind than the (1.32.A) to (1.32.D), namely the

\bigskip

\noindent {\bf Theorem 9.} {\it We have the implication}
$$
S_u \, \widetilde M^3 (G) \in GSC \Longrightarrow G \in QSH.
$$

With this, Theorems 8 and 9 put together clinch the proof of our Theorem 3.

\smallskip

Theorem 8 is the object of Part [II] of the trilogy and we discuss it in this section. Theorem 9 is the object of Part [III] and it will be very briefly discussed in the next section of this survey.

\smallskip

Before discussing holes and ditches, let us return to the context of the bicollared handles, their partial foliations and the induced walls, the spare parts of the $X^2$ in (1.6), the $2^{\rm d}$ REPRESENTATION space of $G$.

\smallskip

We have the compact walls $W ({\rm BLUE})$, $W ({\rm RED})$, $W ({\rm BLACK})$ homeomorphic, respectively to $S^2$, $S^1 \times [0,1]$ and $D^2$ and then we also have the non-compact security walls $W_{\infty}$(BLACK). When we will write $W_{(\infty)}$(BLACK) that will mean ``$W$(BLACK) and/or $W_{\infty}$(BLACK)''. Each $W$(RED) has two parts one contained inside the bicollared $3^{\rm d}$ handles of index $\lambda = 0$, call them $H^0$, and another one outside them, coming with a SPLITTING $W ({\rm RED}) = W ({\rm RED} \cap H^0) \cup W ({\rm RED} - H^0)$. 

\smallskip

The limiting position of the $W$(BLUE), $W$(RED) respectively, are denoted $S_{\infty}^2$(BLUE), $(S^1 \times I)_{\infty}$(RED). With these things we have the following formula, a precise version of (1.11.2),
$$
f \, {\rm LIM} \, M_2 (f) = f \left[ \sum W ({\rm RED}) + W_{(\infty)} ({\rm BLACK}) \right] \cap \left[ \sum S_{\infty}^2 ({\rm BLUE}) \cup \sum (S^1 \times I)_{\infty} ({\rm RED}) \right] , \leqno (1.33)
$$
and more precisely we have $W ({\rm RED}) \cap S_{\infty}^2 ({\rm BLUE}) \subset W ({\rm RED} \cap H^0)$.

\smallskip

There will be two kinds of Holes, the normal holes, very much like in the toy model, and the BLACK Holes which {\bf include} the items accompanying the $p_{\infty\infty}$'s but also more mundane usual Holes, living inside the compact $W ({\rm BLACK})$'s; they are treated on par with the normal holes from formula (1.34) below. The $p_{\infty\infty}$'s live inside the $W_{(\infty)} ({\rm BLACK})$. When we will write ``$-H$'', this will means ``with {\bf all} Holes deleted''. But then, deleting $p_{\infty\infty}$'s means deleting closed subsets while deleting the other Holes (like in (1.34) or the mundane BLACK ones), means deleting open sets. And any deletion is accompanied by addition of compensating 2-handles, like our $D^2 (p_{\infty\infty})$ for the $p_{\infty\infty}$'s, and, generically, these all are our $D^2 (H)$ (like in (1.32.D)). And it should be understood that $\{ D^2 (p_{\infty\infty})\} \subset \{ D^2 (H) \}$.

\smallskip

With this, things, in connection to (1.33) we will have
$$
\{\mbox{Normal Holes}\} \subset \sum W ({\rm BLUE}) + \sum W ({\rm RED} - H^0), \leqno (1.34)
$$
$$
\mbox{and dually}, \ \mbox{DITCHES} \subset \sum W ({\rm RED} \cap H^0) + \sum W_{(\infty)} ({\rm BLACK})  . 
$$

It is the $W_{(\infty)} ({\rm BLACK})$  which carry the BLACK holes, remember. With ditches, Holes and partial ditch-filling, there is now an idea of {\bf high-dimensional geometric realization of the zipping} $X^2 \underset{f}{\longrightarrow} fX^2$ the big quotient space projection from (1.13). This will be now the real life version of the Lemma 6 from the Toy-Model.

\smallskip

The offshot of this geometric realization is another functor $S_b^{(')}$, besides the $S_u^{(')}$. In a first approximation, the high-dimensional realization of the zipping replaces the big quotient space projection by an infinite sequence of inclusion maps, much more manageable then the quotient projections. But other elementary steps, the so-called ditch-jumping steps, brought about by the triple points $M_3(f)$, are necessary too. And here are some important items concerning the infinitely many steps of the geometric, high-dimensional realization of the zipping:
\begin{enumerate}
\item[$\bullet$)] All the elementary steps concerned are GSC-preserving.
\item[$\bullet$$\bullet$)] All the action takes place in the {\bf additional} $N$ dimensions.
\end{enumerate}

\smallskip

So, we start by getting an $(N+4)$-dimensional cell-complex
$$
S_b^{(')} (\widetilde M^3 (G) - H)
$$
with functoriality properties analogous to the (1.32.A), (1.32.B), (1.32.C), and everything is again equivariant. We will not bother to rewrite for 
$S_b^{(')} (\widetilde M^3 (G)-H)$, the formulae like (1.32.A) to (1.32.C).

\smallskip

The $S_b^{(')} (\widetilde M^3 (G) - H)$ comes equipped with PROPER embeddings
$$
\sum_{\mbox{\footnotesize all $H$'s}} C(H) \xrightarrow{ \ \beta \ } S_b^{(')} (\widetilde M^3 (G)-H) 
$$
and now one {\bf defines} (or rather one imposes, by definition) a reconstruction formula

\bigskip

\noindent (1.35) \quad $S_b^{(')} (\widetilde M^3 (G)) \equiv S_b^{(')} (\widetilde M^3 (G) - H) + \underset{H}{\sum} \, D^2 (H)$ (added this time along $\beta$ not along $\alpha$ like in (1.32.D)).

\bigskip

It is the $\beta \, C(H)$ for the $H \in \{$normal Holes (1.34)$\}$, which are involved in the process of PARTIAL DITCH FILLING, and the toy-model Figure 8 should help in understanding what this means.

\smallskip

The mundane BLACK holes do not partake into this {\bf partial} ditch filling. But, notice that, while in (1.32.D) the reconstruction was a fact to be proved, here for $S_b^{(')}$, it is essentially a definition. But, again, the $S_b^{(')}$ analogues of (1.32.A), (1.32.B), (1.32.C) hold for the functor $S_b^{(')}$. Now here comes the key fact: because the construction of $S_b^{(')}$ proceeds by GSC preserving steps (like inclusions which come from dilatations, or additions of handles of index $>1$) and {\bf not} by quotient-space projection, we have now, the following item:

\bigskip

\noindent {\bf Lemma 10.} 1) {\it $S'_b \, \widetilde M^3 (G)$ is GSC.}

\medskip

2) {\it Using a relatively easy transformation, proceeding on the same lines as in {\rm 2)}, Theorem $8$, i.e. via healing,
$$
S'_b \, \widetilde M^3 (G) \Longrightarrow S_b \, \widetilde M^3 (G) \, ,
$$
the $S_b \, \widetilde M^3 (G)$ is also GSC.}

\bigskip

A priori it looks as if the Lemma 10 should be quite automatic. But this is not quite so. In getting the map $\beta$ to be PROPER we encounter difficulties, coming on the one side form $M_3 (f) \ne \emptyset$ and, most importantly, from the possible {\bf bad cycles} created by the intersection
$$
\emptyset \ne (\mbox{collapsing flow lines}) \pitchfork (\mbox{zipping flow lines}) \subset X^2 \, .
$$
And handling these additional problems makes that the proof of Lemma 10 is not trivial.

\smallskip

Here comes now a first connection between our two functors.

\bigskip

\noindent {\bf Lemma 11.} {\it There is an {\ibf simple-minded diffeomorphism} $\eta$ entering in the diagram below:
$$
\xymatrix{
S'_b (\widetilde M^3 (G)-H) \ar[rrr]^-{\eta}_-{\mbox{\footnotesize simple-minded diffeomorphism}} &&&S'_u (\widetilde M^3 (G)-H) \\
&\underset{\mbox{\footnotesize all $H$'s}}{\sum} C(H) . \ar[ul]^{\beta} \ar[urr]_{\alpha}
}
\leqno (1.36)
$$

This diagram {\ibf commutes up to homotopy}.

\smallskip

The ``simple-minded'' here, is to be understood in opposition to the quality which the later diffeomorphism $S_b (\widetilde M^3 (G)) \underset{\rm DIFF}{=} S_u (\widetilde M^3 (G))$ may have. And this one is not simple-minded, by {\ibf any} means.}

\bigskip

It is {\bf not} claimed that this homotopy is PROPER, otherwise we would could get isotopy too, then we would have $S'_b (\widetilde M^3 (G)) \underset{\rm DIFF}{=} S'_u (\widetilde M^3 (G))$ quite cheaply, and of course with this diffeomorphism we would be done, as far as Theorem 8 is concerned.

\smallskip

But life is more complicated than that. So, let us start by going back to (1.36), which has to be compared to the lower triangle in (1.30), of course. Notice also that the links $\alpha \left( \underset{H}{\sum} \, C(H)\right)$, $\beta \left( \underset{H}{\sum} \, C(H)\right)$ occurring in (1.36) come with canonical-framings.

\smallskip

We move now from $\widetilde M^3 (G)$ downstairs to the compact space $M^3 (G)$.

\bigskip

\noindent {\bf The compactness Lemma 12.} {\it When we move from {\rm (1.36)} to the corresponding diagram {\ibf downstairs}
$$
\xymatrix{
S'_b ( M^3 (G)-H) \ar[rrr]^-{\eta}_-{\mbox{\footnotesize simple-minded diffeomorphism}} &&&S'_u ( M^3 (G)-H) \\
&\underset{\mbox{\footnotesize all $H$'s}}{\sum} \{C(H) , \ \mbox{with framing}\} \ar[ul]^{\beta} \ar[urr]_{\alpha}
}
\leqno (1.37)
$$
then this diagram commutes now up to PROPER homotopy.}

\bigskip

It is the compactness of $M^3 (G)$ which makes this possible. But the spaces $S'_b ( M^3 (G)-H)$, $S'_u ( M^3 (G)-H)$ (nor the $S'_b ( M^3 (G))$ or $S'_u ( M^3 (G))$) are certainly NOT compact.

\bigskip

We will discuss this all-important Lemma 12 in the next section, but right now, here is how the Theorem~8 is proved, modulo all the things said.

\smallskip

To begin with, for $S'_b (M^3 (G))$, $S'_u (M^3 (G))$ we have reconstruction formulae like (1.35), (1.32-D), but now downstairs. Next, as a consequence of Lemma 12 we also find now a diffeomorphism, downstairs of course,
$$
S'_b (M^3(G)) \underset{\rm DIFF}{=} S'_u (M^3(G)) \, . \leqno (1.38)
$$
The objects involved here are not manifolds but they make sense as closed subsets (actually subcomplexes) of some high-dimensional euclidean space, so one can happily talk about diffeomorphisms in our singular context.

\smallskip

Of course, also, the implication Lemma 12 $\Rightarrow$ formula (1.38), uses the fact that in the high dimension $N+4$ we also have the implication
$$
\mbox{PROPER homotopy} \Longrightarrow \mbox{isotopy, for links.}
$$

We insist, once more, that the objects involved in (1.37) are not compact.

\smallskip

[The qualification ``COMPACTNESS'' for Lemma 12 stems from the fact that, because $M^3(G)$ {\bf is} compact, we have some {\bf compactifications} for the $(N+4)$-dimensional non-compact cell-complexes, which occur in (1.37). This will be discussed later on.]

\bigskip

Consider now the following diagram which is given to us by FUNCTORIALITY
$$
\downarrow\!\!\!\!^{\overset{\mbox{$u$}}{\underset{\mbox{$\approx$}}{-\!\!-\!\!-\!\!-\!\!-\!\!-\!\!-\!\!-\!\!-\!\!-\!\!-\!\!-\!\!-\!\!-\!\!-\!\!-\!\!-\!\!-\!\!-\!\!-\!\!-\!\!-\!\!-\!\!-\!\!-\!\!-\!\!-\!\!-\!\!-\!\!-\!\!-\!\!-\!\!-\!\!-\!\!-\!\!-\!\!-\!\!-\!\!-\!\!-\!\!-\!\!-\!\!-\!\!-\!\!-\!\!-\!\!-\!\!-\!\!-\!\!-\!\!-\!\!-\!\!-\!\!-\!\!-\!\!-\!\!-\!\!-\!\!-\!\!-\!\!-\!\!-\!\!-\!\!-\!\!-\!\!-\!\!-\!\!-\!\!-\!\!-\!\!-\!\!-\!\!-\!\!-\!\!-\!\!-\!\!-\!\!-\!\!-\!\!-\!\!-\!\!-\!\!-\!\!-\!\!-\!\!-\!\!-\!\!-\!\!-\!\!-\!\!-\!\!-\!\!-\!\!-}}}\!\!\!\!\downarrow
$$
\vglue-7mm
$$
\xymatrix{
S_u (\widetilde M^3(G)) \Longleftarrow S'_u (\widetilde M^3 (G)) \ar[r]^{\mathcal V}_{\approx} \ar[d]^p &S'_b (\widetilde M^3 (G)) \Longrightarrow S_b (\widetilde M^3(G)) \, . \ar[d]^p \\
S'_u (M^3(G)) \ar[r]_{\approx}^{\mathcal W} &S'_b (M^3(G))
}
\leqno (1.39)
$$
Here, the bottom arrow ${\mathcal W}$ is the diffeomorphism (1.38) and the two vertical maps $p$ are the natural universal covering space projections. Hence, because ${\mathcal W}$ is a diffeomorphism, the induced map ${\mathcal V}$ is also a diffeomorphism.

\smallskip

Now, by Lemma 10, $S'_b (\widetilde M^3 (G)) \in {\rm GSC}$ hence $S'_u (\widetilde M^3 (G))$ is also GSC. [Notice that this concerns the upper line of our (1.39), the objects on the lower line come with
$$
\pi_1 \, S'_u (M^3 (G)) = \pi_1 \, S'_u (M^3 (G)) = G \, . ]
$$

Finally, the two double arrows which occur  in the top line of (1.39) are transformations which, as transformations, are isomorphic, irrespective of what their respective sources or targets may be. The claimed isomorphism is then between the transformations, as abstract operation. But then, the diffeomorphism ${\mathcal V}$ induces a diffeomorphism $u$, and invoking again Lemma 10 our $S_b (\widetilde M^3 (G))$ is GSC, hence so is $S_u (\widetilde M^3(G))$, and Theorem 8 is proved.

\bigskip

\noindent {\bf A philosophical comment.} -- An object like $S_u^{(')} (\widetilde M^3 (G))$, and which is certainly non-compact, has two kinds of infinities. There is, to begin with, the big mysterious INFINITY of the group $G$, considered as a metric space. Asymptotic properties of $G$ like QSF, $\pi_1^{\infty}$, a.s.o, really concern this INFINITY.

\smallskip

But then $S_u^{(')} (\widetilde M^3 (G))$ has, in addition to this God-given global INFINITY, a more artificial local, man-made infinity. The $S_u^{(')}(\widetilde M^3 (G))$ consist of local units corresponding to things like the bicollared handles $H_i^{\lambda} (\gamma)$ from Figure 4, actually $(N+4)$-dimensional {\bf infinitely foamy}  $(N+4)$-dimensional thickening of some very dense $2^{\rm d}$-skeleton of the already non-compact $H_i^{\lambda} (\gamma)$'s, hence they are quite violently non-compact. This is, essentially, our artificial man-made infinity.

\smallskip

Now, in going down to $S'_u (M^3 (G) - H)$, $S'_b (M^3 (G) - H)$, we get rid of the God-given INFINITY of $G$, which is hidden in the fibers of the universal covering map $p$ which occurs in (1.39).

\smallskip

But then, when we want to get the diffeomorphism (1.38), meaning when we want to prove the COMPACTNESS Lemma 12, we still have to come to grip with the man-made, lesser infinity. This is what we will discuss next, after some additional explanations concerning the compactness lemma.

\section{Additional details concerning the compactification and Part III of the trilogy}

Corresponding to the two halves of the title of this section 2, we will have now two subsections. The first one is

\subsection{Additional details concerning the compactification}

We will discuss now some of the ideas behind the proof of Lemma 12, and, for this, we move downstairs to $S'_u (M^3 (G)-H)$.

\smallskip

For the various holes $H_1 , H_2, \ldots$ downstairs, now at level $\{ \Theta^3 (fX^2-H)/G$, which includes from now all the refinements like the ones discussed in (1.31) to (1.31.B)$\}$ and also the $\underset{i}{\Sigma} \, C(H_i)$, which occur in diagram (1.37). What we know, a priori, is that the diagram in question commutes up to homotopy and the issue now is to show that this is actually PROPER homotopy. [Keep in mind that all our story is equivariant, at all times. The diagram we talk about now is the equivariant diagram (1.36), upstairs, projected downstairs by quotient space projection.]

\smallskip

We express now the homotopy commutativity of (1.37) as follows.

\bigskip

\noindent (2.1) \quad For each $H_i$ there is an arc $\gamma_i \subset S'_u (M^3 (G)-H)$ joining $\alpha \, C(H_i)$ to $\eta \beta \, C(H_i)$ such that the closed loop
$$
\Lambda (H_i) \equiv \alpha \, C(H_i) \underset{\mbox{\footnotesize $\overbrace{\gamma_i}$}}{\bullet} (\eta\beta \, C(H_i))^{-1}
$$
is null-homotopic in $S'_u (M^3 (G)-H)$. End of (2.1).

\bigskip

This is, of course, text-book level elementary topology, but it is actually not so innocent, since we have now the following

\bigskip

\noindent {\bf Lemma 13.} {\it One can {\ibf choose} the arcs $\gamma_i$ entering the definition of $\Lambda (H_i)$ above, so that the following two things should happen.}

\medskip

1) {\it The lengths $\Vert \Lambda (H_i) \Vert$ are controlled by the zipping length of the REPRESENTATION {\rm (1.6)}, and hence there is a {\ibf uniform bound} $K_1$ which is such that
$$
\Vert \Lambda (H_i) \Vert < K_1 \, , \ \mbox{for all $i$'s}. \leqno (2.2)
$$
}

2) {\it We have
$$
{\lim}_{i \, \to \, \infty} \, \Lambda (H_i) = \infty \, , \leqno (2.3)
$$
in the non-compact space $S'_u (M^3 (G)-H))$.}

\bigskip

All this story, with (2.1), (2.2), (2.3), is valid both downstairs and upstairs, and if we stated it downstairs it is because that is the context where we will make use of it.

\smallskip

And now we turn for good to the downstairs context, where the ambient space $M^3 (G)$ is {\bf compact}. And, using Lemma 13 we can prove the

\bigskip

\noindent {\bf Lemma 14.} 1) {\it There is a way to glue to the infinity of $S'_u (M^3 (G)-H)$ various lines and points living themselves at the infinity of $\Sigma_0 (\infty)$ {\rm (1.31.B)} so as to get a space $S'_u (M^3(G)-H)^{\wedge} \supset S'_u (M^3 (G)-H)$, which has the following property. When we consider the natural embedding
$$
S'_u (M^3(G)-H) \subset \{\mbox{The compact metric space} \ M^3(G) \times B^L , \ \mbox{for some high $L$}\} ,
$$
then $(S'_u (M^3(G)-G))^{\wedge}$ is exactly the closure
$$
\overline{S'_u (M^3(G)-H)} \subset M^3(G) \times B^L \, .
$$
}

2) {\it So $S'_u (M^3 (G)-H)^{\wedge}$ is COMPACT. This is our compactification of $S'_u (M^3 (G)-H)$.}

\medskip

3) {\it For any subsequence
$$
\{ \Lambda_{i_1} , \Lambda_{i_2} , \ldots \} \subset \{\Lambda_1 , \Lambda_2 , \ldots \} \, ,
$$
where $\Lambda_i \equiv \Lambda (H_i)$, we can find a sub-sub-sequence
$$
\{\Lambda_{j_1} , \Lambda_{j_2} , \ldots \} \subset \{\Lambda_{i_1} , \Lambda_{i_2} , \ldots \}
$$
such that there is a closed curve $\Lambda_{\infty} \subset (S'_u (M^3 (G)-H))^{\wedge} - (S'_u (M^3 (G)-H))$, with the property that
$$
{\lim}_{n \, \to \, \infty} \Lambda_{j_n} = \Lambda_{\infty} \, , \ \mbox{{\ibf uniform convergence} in} \ (S'_u (M^3 (G)-H))^{\wedge} \, . \leqno (2.4)
$$
}

\bigskip

The $S'_u (M^3 (G) - H)^{\wedge}$ is a compact metric space and in (2.4) we talk about uniform convergence of sequences of continuous maps from $S^1$ to this metric space.

\smallskip

Once the $S'_u (M^3 (G)-H)$ includes the ingredient (1.31.B) we can put to use the Lemma 14 so as to get the following MAIN LEMMA below. And our step (1.31.B) is certainly necessary here.

\bigskip

\noindent {\bf MAIN LEMMA 15.} {\it For each closed curve $\Lambda_n = \Lambda (H_n)$, see here {\rm (2.4)}, there is a singular disc $D_n^2$ cobounding $\Lambda_n$, $D_n^2 \to S'_u (M^3 (G)-H)$, with $\partial D_n^2 = \Lambda_n$, having the property that
$$
{\lim}_{n \, \to \, \infty} D_n^2 = \infty \quad \mbox{in} \quad S'_u (M^3 (G)-H) \, .
$$
}

\bigskip

\noindent {\bf A remark concerning Lemma 15.} The statement above has a  tint of $\pi_1^{\infty} = 0$ flavour, but no simple connectivity at infinity of any locally compact space is actually proved. $\Box$

\bigskip

From the Main Lemma 15, the compactness Lemma 13 follows, and with this we have finished our discussion of Theorem 8 and we will turn to Theorem 9; we will briefly describe Part III of the Trilogy.

\smallskip

So, here comes now our next subsection.

\subsection{A glimpse into Part III of the trilogy}

\smallskip

At this point let us compare the stabilization lemma from Part I of this survey, to the following classical result of J. Stallings \cite{42}

\bigskip

\noindent (2.5) \quad (J. Stallings) Let $M^n$ be a smooth open contractible manifold. Then, if $N$ is high enough
$$
M^n \times R^N \underset{\rm DIFF}{=} R^{n+N} \, .
$$
The point is that in (2.5), when the multiplying factor is the open $R^N$, $M^n$ could be essentially anything, it could actually be the Whitehead manifold ${\rm Wh}^3$. This is certainly not the case for the {\bf stabilization lemma} $(M^n \times B^p \in {\rm GSC} \Rightarrow M^n$ Dehn-exhaustible), from Part I of our survey. There, unlike in (2.5), the multiplying factor $B^p$ is compact. So, in the context of the stabilization lemma things are {\bf transversally compact}, while in the case of (2.5) this is certainly not the case. And now, for Theorem 9 which is the object of our present discussion, we will need something somehow similar to the stabilization lemma, and the problem is to stay transversally compact; we will call this concern the {\bf Stallings barrier}. And the point is that the non-metrizability barrier already discussed and the Stallings barrier somehow play against each other, but we have to stay on the good side of both of them, at the same time.

\smallskip

So, here is a glimpse of how the proof of Theorem 9 goes. And, at this point we have to say a few more things concerning the $S_u (\widetilde M^3(G))$, a very detailed knowledge of which is actually necessary for the proof of Theorem 9; but only a minimum necessary will be developed here.

$$
\includegraphics[width=120mm]{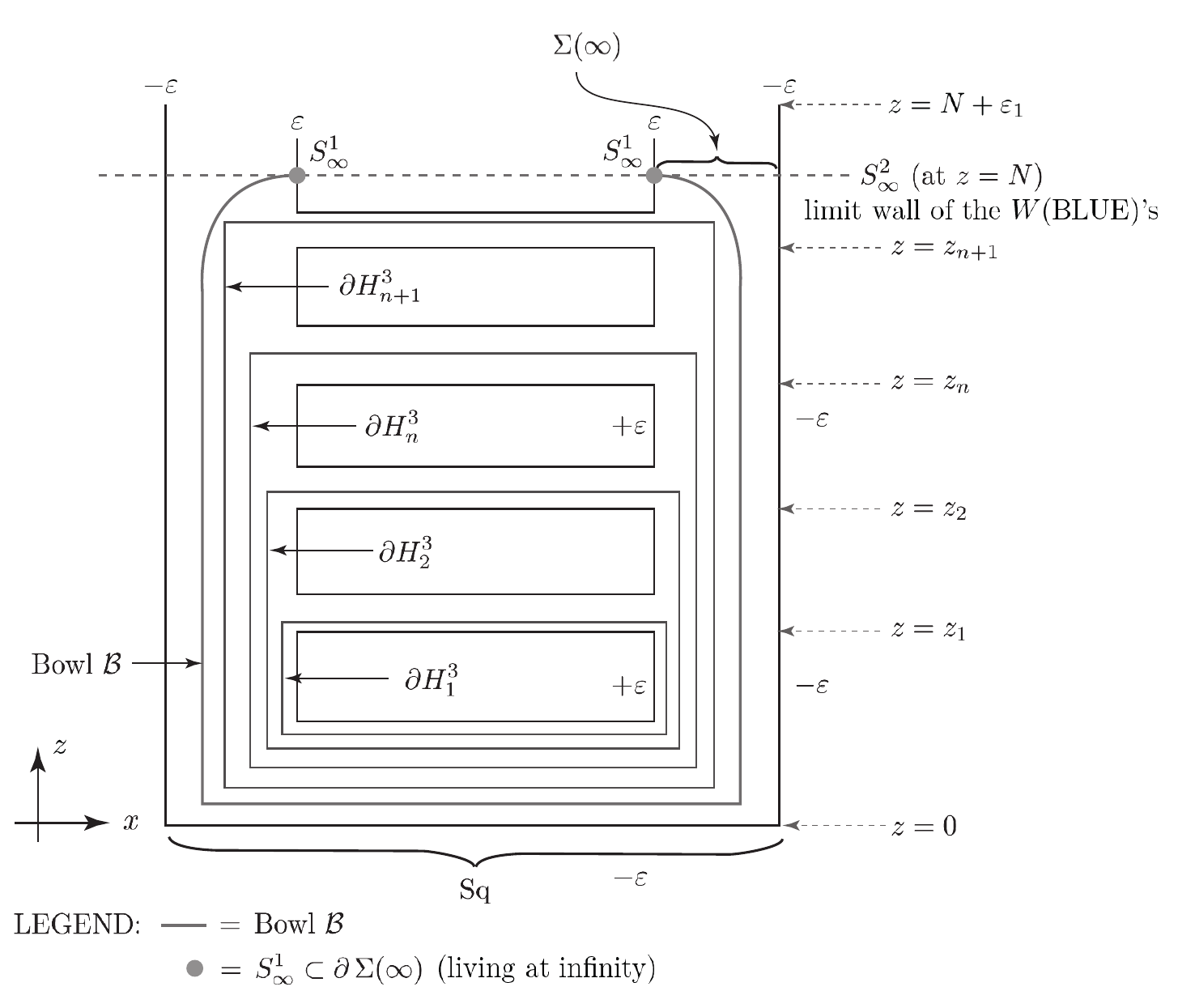}
$$
\label{fig9}
\centerline {\bf Figure 9.} 
\begin{quote}
The $U^3 ({\rm Blue})$ embellished with the BOWL ${\mathcal B}$ and with infinitely many 3-handle attaching spheres $\partial H_1^3 , \partial H_2^3 , \ldots$ which live inside $U^3 ({\rm BLUE}) \cup {\mathcal B} \times [0,\infty)$ and accumulate on $({\mathcal B} \times \{\infty \}) \cup S_{\infty}^2$.
\end{quote}

\smallskip

So, there is, to begin with, a singular, locally finite $3^{\rm d}$ complex $\Theta^3 (fX^2)$. This is not quite the naive $3^{\rm d}$ thickening of the non-locally  finite $fX^2$. Some deletions like in the Figures 5 and 6 are certainly necessary for having a locally finite $\Theta^3 (fX^2)$. For technical reasons, some other embellishments are necessary too. Some are discussed in connection to (1.31) to (1.31.B).

\smallskip

The singularities of $\Theta^3 (fX^2)$ are, essentially, of the undrawable type for which the appropriate figures are provided in \cite{27}, \cite{28}, but see also Part I of this survey, (OR the related cylinders $C(p_{\infty\infty}) \times \left[-\frac\varepsilon4 , \frac\varepsilon4 \right]$, in Figure 6). At this point some of the very initial and elementary technology from \cite{27}, \cite{28} is necessary. For the undrawable singularities of $\Theta^3 (fX^2)$, some kind of very specific {\bf resolution of singularities} (or desingularisation) are available. Like in the classical case, we have here
$$
\left.
\begin{matrix}
\{3^{\rm d} \ \mbox{space without undrawable singularities}\} \\ 
\twoheaddownarrow \\ 
\Theta^3 (fX^2)
\end{matrix}
\right\} \mbox{we call, generically, this desingularization ${\mathcal R}$},
\leqno (2.6)
$$
where the surjection above exactly blows up the undrawable singularities. There are exactly
$$
2^{{\rm card} \, \{\mbox{\footnotesize undrawable singularities of} \ \Theta^3 (fX^2)\}}
$$
desingularizations ${\mathcal R}$ which are possible. And, to each of them a precise $4^{\rm d}$ thickening $\Theta^4 (\Theta^3 (fX^2) , {\mathcal R})$ is prescribed in \cite{27}, \cite{28}. The $\Theta^4 (\Theta^3 (fX^2) , {\mathcal R})$ is still not smooth, nothing is done concerning the various kind of $p_{\infty\infty}$'s and their compensating 2-handles, which are left with singular attachements, and also it is certainly ${\mathcal R}$-dependent. But the point is that when one goes to $\Theta^4 (\Theta^3 (fX^2) , {\mathcal R}) \times B^N$ then the ${\mathcal R}$-dependence gets washed away and (at least in a first approximation) we have actually, as already said
$$
S_u (\widetilde M^3 (G)) = \Theta^4 (\Theta^3 (fX^2) , {\mathcal R}) \times B^N \, ; \leqno (2.7)
$$
and for the present section we will ignore further subtelties and take the (2.7) as our definition of $S_u (\widetilde M^3(G))$. But keep in mind that this is just a simplified version, for expository purposes. The $\Theta^3 (fX^2)$ is certainly $G$-dependent and also $G$-equivariant
$$
G \times \Theta^3 (fX^2) \longrightarrow \Theta^3 (fX^2) \, .
$$
Because of its ${\mathcal R}$-dependence, $\Theta^4 (\Theta^3 (fX^2) , {\mathcal R})$ fails to be equivariant; no group action exists there.

\smallskip

Now, at the level of (2.7) the ${\mathcal R}$-dependence has gone, things are now {\bf canonical}, and so at this level we have again our natural free action
$$
G \times S_u (\widetilde M^3 (G)) \longrightarrow S_u (\widetilde M^3 (G)) \, .
$$

We take now a closer look at $fX^2$. We will describe a typical local piece of $fX^2$ in a simplest possible case and, at the target $\widetilde M^3 (G)$, inside some $f$ (bicollared handle of index $\lambda = 0$, $H_i^0 (\gamma)$). But for simplicity's sake we will describe it as living inside $R^3 = (x,y,z)$. Our local piece is then
$$
U^2 ({\rm BLUE}) = \left\{ \underbrace{[-1 \leq x \leq 1 , -1 \leq y \leq 1 , z=0]}_{\mbox{\footnotesize call this ${\rm Sq}$ like ``square''}} \ \cup \ \left(\partial \, {\rm Sq} \, \times [0 \leq z \leq N+\varepsilon_1] \right) \right\} \ \cup \leqno (2.8)
$$
$\cup \ \{$infinitely many $2$-handles parallel to Sq and BLUE, like it, call them $\underset{i \, = \, 1}{\overset{\infty}{\sum}} \, {\rm Sq} \times \{ z_i \}$, where $0 < z_1 < z_2 < \ldots < N$ and $\underset{n \, = \, \infty}{\lim} \, z_n = N\}$ (see Figure 9).

\bigskip

The mythical $\delta H_i^0 (\gamma)$ (1.8.1) which we call $S_{\infty}^2$, lives at $z=N$ and $U^2 ({\rm BLUE})$ communicates with the rest of $fX^2$ via its outer surface, far from the attachments $\partial \, {\rm Sq} \times \{z_n\}$ of the 2-handles ${\rm Sq} \times \{ z_n\}$. When we move from $fX^2$ to $\Theta^3 (fX^2)$ then, to begin with, the $U^2 ({\rm BLUE})$ is replaced by
$$
U^3 ({\rm BLUE}) = U^2 ({\rm BLUE}) \times [-\varepsilon , \varepsilon] - \{ \partial \, {\rm Sq} \times [(z=N) \times \varepsilon] \} \, , \leqno (2.9)
$$
which occurs as $S_{\infty}^1$ in Figure 9, representing $U^3 ({\rm BLUE})$.

\smallskip

Now, what we see, in a first instance in Figure 9 is an {\bf infinite foaminess} of $U^3 ({\rm BLUE}) \subset \Theta^3 (fX^2)$ and, as a first step towards the proof of Theorem 9, we want to get rid (of most) of it.

\smallskip

We will start by embellishing $U^3 ({\rm BLUE})$ with a PROPER hypersurface

\bigskip

\noindent (2.10) \quad ${\mathcal B}$ (like ``BOWL'') $= \{$a copy of $R^2$, PROPERLY embedded in $U^3 ({\rm BLUE})$ far from $\partial U^3 ({\rm BLUE})$ and resting, at infinity, on $S_{\infty}^1$ (which is not part of $U^3 ({\rm BLUE})$, remember)$\}$ (see Figure 9).

\bigskip

Next, we enlarge $U^3 ({\rm BLUE})$ by going to
$$
U^3 ({\rm BLUE}) \underset{\overbrace{\mbox{\footnotesize ${\mathcal B} \times \{0\} = {\mathcal B}$}}}{\cup} {\mathcal B} \times [0,\infty) \, , \leqno (2.11)
$$
which will be also included into the definition of $\Theta^3 (fX^2)$, on par with (1.31.B).

\smallskip

To (2.11) we add now infinitely many handles of index $\lambda = 3$, in cancelling position with the $\underset{n \, = \, 1}{\overset{\infty}{\sum}} \, {\rm Sq} \times \{z_n\}$, call them
$$
H_1^3 , H_2^3 , H_3^3 , \ldots \, . \leqno (2.12)
$$

Their attaching zones are the $\left\{\underset{n \, = \, 1}{\overset{\infty}{\sum}} \, \partial H_n^3 \ \mbox{suggested in Figure 9}\right\}$, which we let afterwards to climb, partially, on ${\mathcal B} \times [0,\infty)$, so as to achieve the following condition
$$
{\lim}_{n \, \to \, \infty} \partial H_n^3 \subset ({\mathcal B} \times \{\infty\}) \cup S_{\infty}^2 \, . \leqno (2.13)
$$
[The partial climb of the $\underset{n \, = \, 1}{\overset{\infty}{\sum}} \, \partial H_n^3$ on ${\mathcal B} \times [0,\infty)$ is to be performed simultaneously for all $\partial H_n^3$'s. To vizualize this bring, in Figure 9, the Bowl ${\mathcal B}$ much closer to the $\varepsilon$-boundary of $U^3$.]

\smallskip

With this the attachment of (2.12) is PROPER.

\smallskip

Finally, we subject our $U^3 ({\rm BLUE})$ to the transformation, gotten by handle-cancellation
$$
U^3 = U^3 ({\rm BLUE}) \Longrightarrow U^3 \cup {\mathcal B} \times [0,\infty) + \sum_{n \, = \, 1}^{\infty} H_n^3 \Longrightarrow \Bigl\{ U^3 \cup {\mathcal B} \times [0,\infty) \ \mbox{with} \leqno (2.14)
$$
$$
\mbox{all the 2-handles} \ \sum_{n \, = \, 1}^{\infty}  {\rm Sq} \times \{z_n\} \ \mbox{\bf deleted} \Bigl\} \, \equiv U^3 ({\rm new}) \, .
$$

$$
\includegraphics[width=140mm]{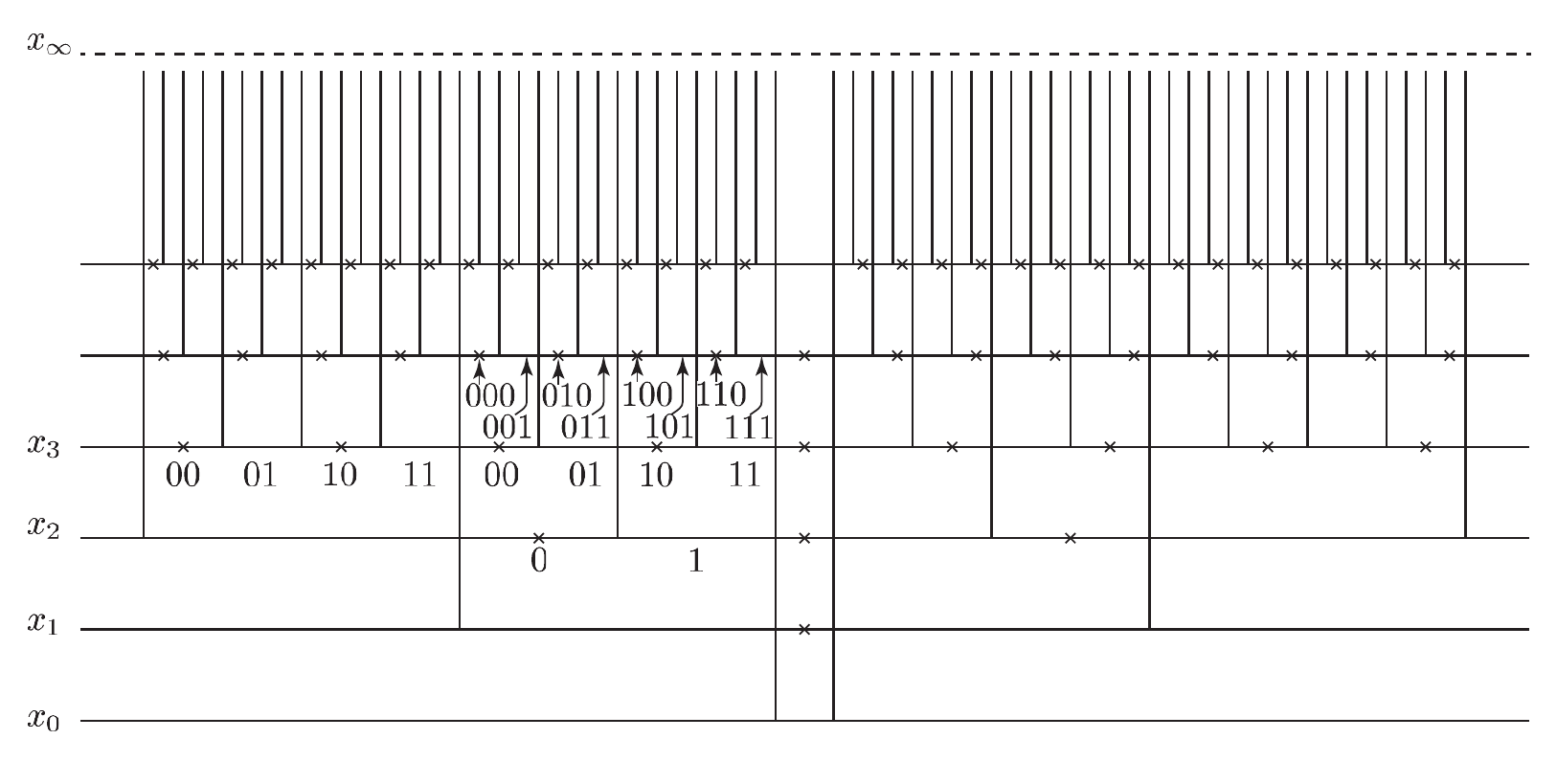}
$$
\label{fig10}
\centerline {\bf Figure 10.} 

\bigskip

When the transformation (2.14) is extended from the $2^{\rm d}$ skeleton in Figure 9 to the whole of the $fX^2 \cap H_i^0 (\gamma)$ (bicollared handle of index $\lambda = 0$), then that 2-skeleton becomes something more complex, like Figure~10 might suggest. But the (2.14) is only a BLUE {\bf game} of type (2.14) played for one bicollared handle of index $0$.

\bigskip

\noindent {\bf Lemma 16.} 1) {\it There is a whole infinite symphony of games like {\rm (2.14)} of all three colours BLUE, RED, BLACK, and taking into account all the infinitely many
$$
fX^2 \cap f (\mbox{bicollared handle} \ H_i^{\lambda} (\gamma)) \subset X^3 \, , \leqno (2.14.1)
$$
and we refer here to the basic REPRESENTATION
$$
X^3 \xrightarrow{ \ f \ } M^3(G) \leqno (2.14.2)
$$
which is locally-finite, equivariant, and of uniformly bounded zipping length.

\smallskip

In {\rm (2.14.1)}, the first $f$ refers to {\rm (1.6)} and $f H_i^{\lambda} (\gamma)$ to {\rm (2.14.2)}. This infinite symphony of games induces a big transformation
$$
\Theta^3 (fX^2) (\mbox{original}) \Longrightarrow \Theta^3 (fX^2) (\mbox{new}) \leqno (2.15)
$$
with the following features.}

\medskip

2) {\it The} (2.15) {\it is $G$-equivariant and, if in {\rm (2.7)} we substitute $\Theta^3 (fX^2) (\mbox{new})$ for $\Theta^3 (fX^2)$, then the $S_u (\widetilde M^3 (G))$ stays GSC.}

\medskip

3) {\it The transformation {\rm (2.17)} stays far from $\Sigma (\infty)$ introduced immediately after {\rm (1.31)} and occurring also in Figure {\rm 5}, and also far from the $\Bigl( p_{\infty\infty} , D^2 (p_{\infty\infty})$ $\times$ $\left[ -\frac\varepsilon4 ,\frac\varepsilon4 \right] \Bigl)$ (see Figure $6$).}

\medskip

4) {\it What the transformation {\rm (2.15)} achieves, is to remove, or obliterate the infinite foaminess of $\Theta^3 (fX^2)$ (original), and we mean here the kind of infinite foaminess which a figure like $10$ may suggest. This feature is important for the next step which follows.}

\bigskip

\noindent [{\bf Technical Remark.} The infinite foaminess is still destroyed if, let us say in the context of a figure like 9 we leave a few (certainly {\bf finitely} many) ${\rm Sq} \times \{z_n \}$'s alive. For technical reasons, we will have to do that, in the context of (2.15) without contradicting equivariance or anything else. In Figure 11, such a finite collection of residual walls is suggested. End of Remark.]

\bigskip

From now on, the $S_u (\widetilde M^3 (G))$ will be conceived on the lines of 2) in the lemma above, i.e. as
$$
S_u (\widetilde M^3 (G)) = \Theta^4 (\Theta^3 (fX^2) ({\rm new}) , {\mathcal R}) \times B^N \, . \leqno (2.16)
$$

$$
\includegraphics[width=130mm]{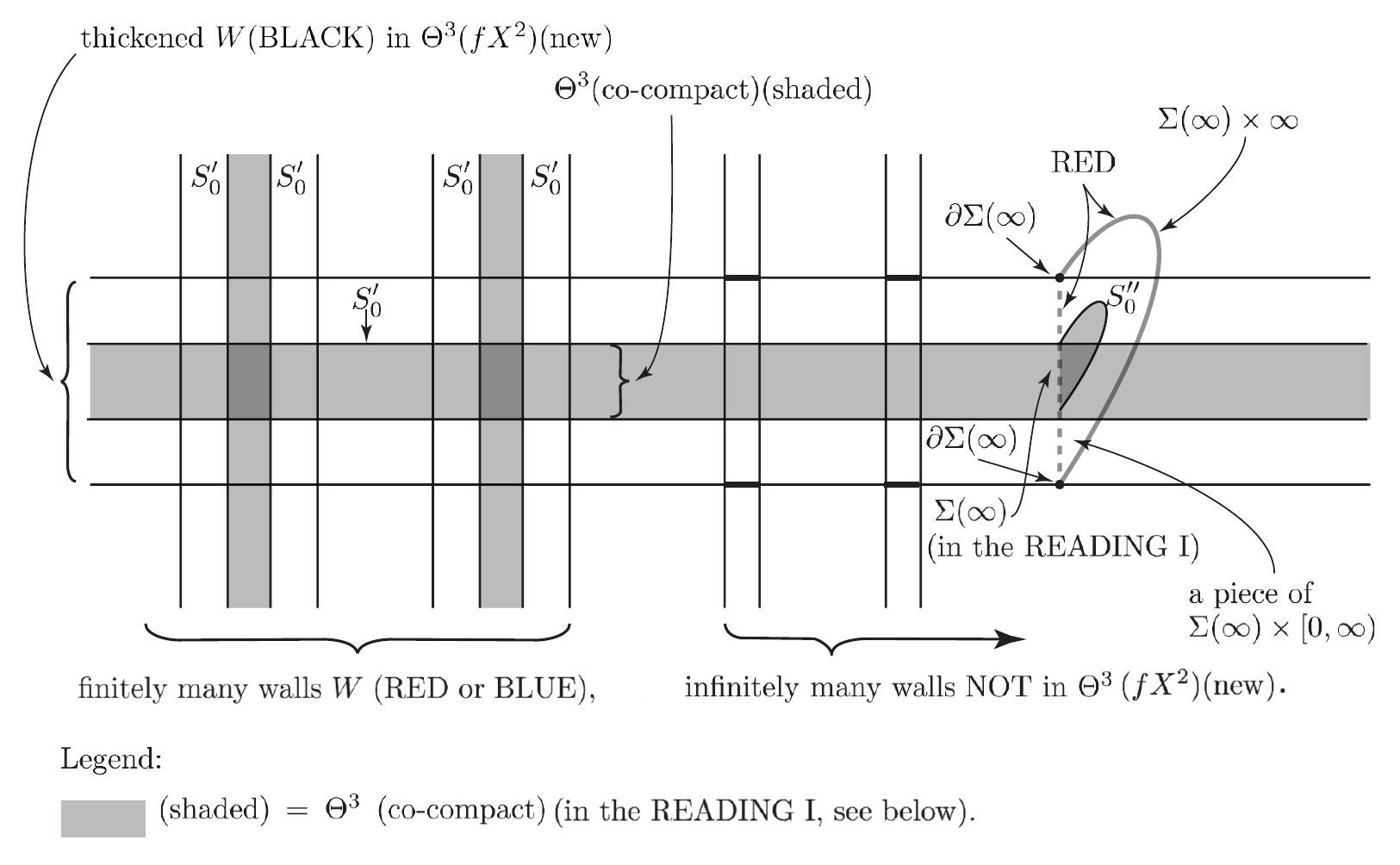}
$$
\label{fig11}
\centerline {{\bf Figure 11.} We are here far from $\{p_{\infty\infty}\}$.}

\bigskip

\noindent {\bf Lemma 17.} 1) {\it Inside $\Theta^3 (fX^2) (\mbox{new})$ lives a non-compact subcomplex
$$
\Theta^3 (\mbox{co-compact}) \subset \Theta^3 (fX^2) (\mbox{new}) \, ,
$$
which is $G$-equivariant and which stays far from those things living close to the infinity of $\Theta^3 (fX^2) (\mbox{new})$ and which prevent the action $G \times \Theta^3 (fX^2) (\mbox{new}) \to \Theta^3 (fX^2) (\mbox{new})$ from being co-compact. So, the action
$$
G \times \Theta^3 (\mbox{co-compact}) \to \Theta^3 (\mbox{co-compact}) \leqno (2.17)
$$
{\ibf is} co-compact, i.e. it comes with a {\ibf compact fundamental domain}.}

\medskip

2) {\it There is an equivariant collapse
$$
\Theta^3 (fX^2) (\mbox{new}) \overset{\pi}{\underset{\mbox{\footnotesize COLLAPSE}}{-\!\!\!-\!\!\!-\!\!\!-\!\!\!-\!\!\!-\!\!\!-\!\!\!-\!\!\!-\!\!\!-\!\!\!-\!\!\!-\!\!\!\longrightarrow}} \Theta^3 (\mbox{co-compact}) \, . \leqno (2.18)
$$
}

\bigskip

The proof is rather technical, but Figure 11 should give an intuitive idea of what $\Theta^3 (\mbox{co-compact})$ looks like. There are two ways to read this figure. In Reading I, the plane of Figure 11 is transversal both to the $W$(BLACK) and to $\Sigma (\infty)$. In this Reading I, $\Sigma (\infty)$ is the ``RED'' detail on the RHS of the figure. In the Reading II what we see in Figure 11 is a piece of $\Sigma (\infty)$ itself. In that reading, the shaded part is then $\Sigma (\infty) \cap \Theta^3 (\mbox{co-compact})$.

\smallskip

Anyway, $\Theta^3$ (co-compact) stops at $S'_0$, $S''_0$ in Figure 11. Actually, the little detail containing $S''_0$ and sticking out of the rest, is a further embellishment for the real-life $(\Theta^3)^{(')}$. This detail consists of some small $3^{\rm d}$ dilatations, the ``fins'', made necessary by the interaction of $\Sigma (\infty)$ with $M_3 (f)$. But we stayed away of these technicalities in this survey. $\Box$

\bigskip

We have already mentioned that the {\bf stabilization lemma}, which was stated and proved in Part I of this survey extends to the case of pure cell-complexes, and so does the property of Dehn-exhaustibility, of course. It is that context which we will use now. With this, we have the following

\bigskip

\noindent {\bf Lemma 18.} {\it The pure cell-complex $\Theta^4 (\Theta^3 (fX^2) (\mbox{new}) , {\mathcal R})$ which occurs in {\rm (2.16)} is Dehn-exhaustible.}

\bigskip

This follows immediately by putting together Theorem 8 from which we get that
$$
S_u (\widetilde M^3 (G)) = \Theta^4 (\Theta^3 (fX^2) ({\rm new}) , {\mathcal R}) \times B^N \in {\rm GSC}
$$
and the extended {\bf stabilization lemma}, which tells us that
$$
\mbox{GSC in the stable range} \Longrightarrow \mbox{Dehn exhaustibility.}
$$

\noindent [A word of caution here concerning the extension of the $n=3$ version of the stabilization lemma, explicitly proved in \cite{26} (for this $n=3$ case) to the stabilization lemma in its high-dimensional incarnation used now, and to its further extension from manifolds to cell-complexes. In the context of \cite{26}, when it was the case of rendering something like the map $Z_{\infty}^n \xrightarrow{ \ \pi_{\infty} \, = \, \pi \mid Z_{\infty}^n \ } M^n$, used in Part I of the present survey when the stabilization lemma was proved, simplicial non-degenerate, then in \cite{26} some jiggling with local elementary riemannian geometry was used. It was even very handy, in the $3^{\rm d}$ context of \cite{26} to invoke affine structures, gotten by appealing to the Smale-Hirsch theory. But analogous elementary arguments can be easily put up in a context like $\pi_{\infty}$, both for our stabilization lemma and for its extension.]

\smallskip

In a nutshell now, here is our strategy for proving the Theorem 9. It consists in the following chain of implications $S_u (\widetilde M^3 (G)) \in {\rm GSC} \Rightarrow \Theta^4 (\Theta^3 (fX^2) ({\rm new}) , {\mathcal R}) \in$ Dehn exhaustible (which we write DE) $\Rightarrow \Theta^3 (fX^2) ({\rm new}) \in {\rm DE} \Rightarrow \Theta^3 (\mbox{co-compact}) \in {\rm QSF} \Rightarrow G \in {\rm QSF}$.

\bigskip

\noindent {\bf Lemma 19.} {\it Moving from dimension $4$ to dimension $3$, we have an implication
$$
\Theta^4 (\Theta^3 (fX^2) (\mbox{new}) , {\mathcal R}) \ \mbox{is Dehn-exhaustible} \ \Longrightarrow \Theta^3 (fX^2) (\mbox{new}) \ \mbox{is Dehn-exhaustible}. \leqno (2.19)
$$
}

\bigskip

We will only give a hint of the proof. Our $\Theta^3 (fX^2) ({\rm new})$ has both undrawable singularities like in \cite{27}, \cite{28} and also singular attachments of 2-handles like the $D^2 (p_{\infty\infty}) \times \left[ -\frac\varepsilon4 ,\frac\varepsilon4 \right]$ in Figure 6. This induces a splitting
$$
\Theta^3 (fX^2) ({\rm new}) = [\Theta^3] \underset{X^2}{\cup} \sum D^2 (p_{\infty\infty}) \times \left[ -\frac\varepsilon4 ,\frac\varepsilon4 \right] \, , \leqno (2.20)
$$
with all the undrawable singularities concentrated in $[\Theta^3]$. Then (2.20) induces a $4^{\rm d}$ splitting
$$
\Theta^4 (\Theta^3 (fX^2) ({\rm new}) , {\mathcal R}) = Y^4 \underset{X^3}{\cup} Z^4 \, , \leqno (2.21)
$$
suggested in Figure 12. We can also see there a foliation ${\mathcal F}$, with $3^{\rm d}$ leaves, which connects $\Theta^4 (\Theta^3 (fX^2) ({\rm new}) , {\mathcal R})$ to $\Theta^3 (fX^2) ({\rm new})$.

\smallskip

On uses then this foliation and the retraction of $\Theta^4 (\Theta^3 (fX^2) ({\rm new}) , {\mathcal R})$, on $\Theta^3 (fX^2) ({\rm new})$ in an argument like in the proof of the stabilization lemma, but a bit more sophisticated, in order to prove what we want. $\Box$

\bigskip

$$
\includegraphics[width=120mm]{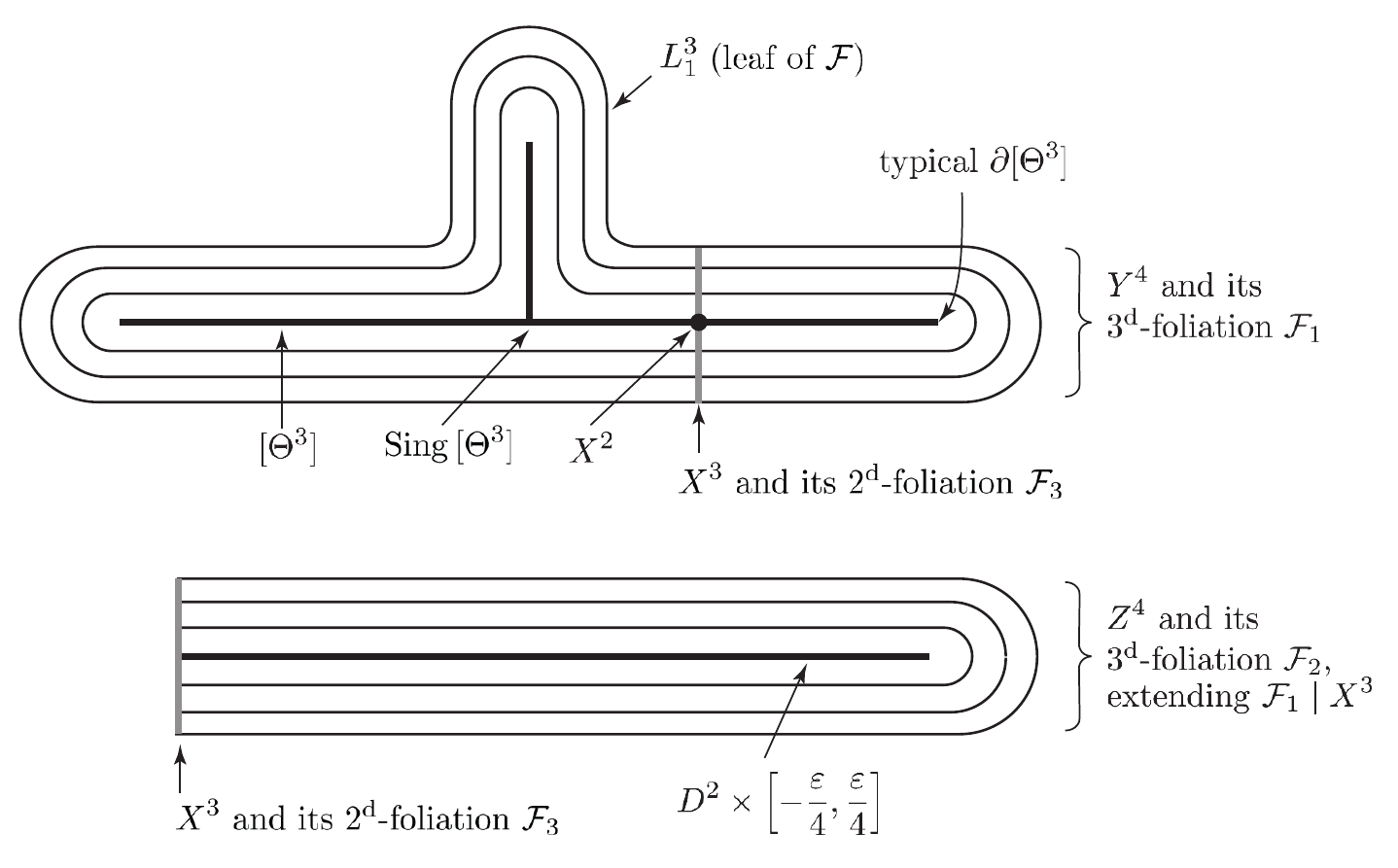}
$$
\label{fig12}
\centerline {{\bf Figure 12.}}
\begin{quote}The splitting loci $X^2 , X^3$ are not our REPRESENTATION spaces. The two halves are to be put together along the common $X^3$.
\end{quote}

\bigskip

\noindent {\bf Lemma 20.} {\it We also have the implication
$$
\Theta^3 (fX^2) (\mbox{new}) \ \mbox{is Dehn-exhaustible} \ \Longrightarrow \Theta^3 (\mbox{co-compact}) \in QSF. \leqno (2.22)
$$
}

\bigskip

Not much details concerning the proof can be given here; and certainly we have not really explained how the $\Theta^3 (\mbox{co-compact})$ is getting generated either. But here are still some ideas concerning our proof. We need to show that for every compact $k \subset \Theta^3 (\mbox{co-compact})$ there is a commutative diagram
$$
\xymatrix{
k \ar[rr]^{j_0} \ar[rd]_-i &&K_0^3 \ar[ld]^-{f_0} \, , \\
&\Theta^3 (\mbox{co-compact})
}
\leqno (2.23)
$$
where $K_0^3$ is compact simply-connected, and where $f_0$ is a simplicial map verifying the Dehn condition $j_0 \, k \cap M_2 (f_0) = \emptyset$.

\smallskip

Here $k$ is given and, what we know is that there is also another commutative diagram
$$
\xymatrix{
\mbox{\hglue -25mm} \qquad \qquad \qquad \quad \qquad \qquad \qquad \qquad \qquad \qquad k \overset{i}{\subset} \Theta^3 (\mbox{co-compact}) \ar[rd]_-j & \subset &\Theta^3 (fX^2)({\rm new}) \\
&K^3 \ar[ru]_-f
}
\leqno (2.24)
$$
where $f$ is now an immersion satisfying the Dehn condition
$$
jk \cap M_2 (f) = \emptyset \, .
$$
The collapse (2.18) is important at this point and, when one looks at things in more detail then one knows that there is a surface with branchings (non-manifold points) properly and PROPERLY embedded
$$
(S , \partial S) \subset (\Theta^3 (fX^2)({\rm new}) , \ \partial \, \Theta^3 (fX^2)({\rm new})) \, , \leqno (2.25)
$$
SPLITTING things away from $\Theta^3 (\mbox{co-compact}) \subset \Theta^3 (fX^2)({\rm new})$. 

\smallskip

Our constant convention for terminology is that ``proper'' means boundary to boundary and interior to interior, while ``PROPER'' means the usual inverse image of compact is compact. Let us call $\Theta^3_0$ these things split away, and which are not in $\Theta^3 (\mbox{co-compact})$. We have, in the context of (2.18), the collapse, restriction of (2.18),
$$
\Theta^3_0 \xrightarrow{ \ \pi \ } S \, , \leqno (2.26)
$$
There is, also, a collection of properly embedded trees
$$
\sum_1^{\infty} A_i \underset{\mbox{\footnotesize PROPER} \atop \mbox{\footnotesize embedding}}{\subset} S \, .
$$
breaking $S$ into manageable compact collapsible pieces which we generically call $B_j$. The $\underset{1}{\overset{\infty}{\sum}} \, \pi^{-1} \, A_i \subset \Theta^3_0$ splits $\Theta^3_0$ into the compact collapsible pieces $\pi^{-1} \, B_j$. In order to go from the given (2.24) to the desired (2.23) involves, among other things, getting from the $K^3$ (2.24) to a $K_0^3 \cap \Theta^3_0 = \emptyset$. This involves finitely many steps, of which we will only describe the simplest. Let us say we find some non-void intersection $K^3 \cap (\pi^{-1} \, A_i - A_i) \ne \emptyset$, which we want to demolish. We may always assume $k$ to be connected and so, in the splitting of $K^3$ by $\pi^{-1} \, A_i - A_i$ there is a piece $jk \subset K_{\alpha}^3$ and a piece $K_{\beta}^3 \cap jk = \emptyset$.

\smallskip

Notice, anyway, that $k \cap \Theta_0^3 = \emptyset$. As a step on the road from (2.24) to (2.23), we will replace $K^3$ by the middle object $K_{\alpha}^3 \cup \{$the cone $\ldots\}$ occuring in the diagram below
$$
k \subset K_{\alpha}^3 \cup \{\mbox{the cone in $B_i$ over $K^3 \cap (\pi^{-1} \, A_i - A_i)$}\} \xrightarrow{ \ \ \ f_{\alpha} \ \ \ } \Theta^3 (fX^2)({\rm new}) \, .
$$
Here $\pi_1 (K_{\alpha}^3 \cup \{\mbox{the cone} \ \ldots \}) = 0$, and $f_{\alpha}$ is a {\bf simplicial} map which satisfies the Dehn condition, but is no longer immersive.

\smallskip

So, in the process, we lost the classical Dehn exhaustibility and we only get the coarser QSF condition (which is Dehn for the now simplicial $f_{\alpha}$ which is no longer immersive). But this brings us also one step closer from (2.24) towards the desired (2.23). The $f_{\alpha}$ is a bit more inside $\Theta^3 (\mbox{co-compact})$ then $f$ was.

\smallskip

This is the general idea.

\subsection{The road-map}

For the convenience of the reader we present here a road-map for our whole argument.

\medskip

1) For a finitely presented group $G$ we consider presentations. And, for our purposes, what we need now are not the well-known 2-complexes built from generators $+$ relations, but some compact 3-manifolds $M^3 (G)$ which have to be {\bf singular} (unless $G$ itself is the fundamental group of a smooth compact 3-manifold). In [I] it is explained why our presentation of $G$ has to be, exactly, 3-dimensional.

\medskip

2) Then one has REPRESENTATIONS of $G$, which are non-degenerate simplicial maps of some ``representation space'' $X$ of dimension 2 or 3, $X \xrightarrow{ \ f \ }\widetilde M^3 (G)$, maps which can be ZIPPED. This is described in Part I of our survey. In that same part we described the $3^{\rm d}$ REPRESENTATION theorem (Theorem 4, Part I), stating that for any group $G$ we can get a REPRESENTATION $X^3 \xrightarrow{ \ f \ }\widetilde M^3 (G)$ which is equivariant, with a locally finite $X^3$ and with a zipping length which is uniformly bounded. Note here that a generic REPRESENTATION space is neither locally finite, nor is the REPRESENTATION equivariant (actually no $G$-action on $X$). So we move now to the present Part II.

\medskip

3) Using 2) one also gets a $2^{\rm d}$ REPRESENTATION
$$
X^2 \xrightarrow{ \ f \ } \widetilde M^3 (G)
$$
with all the three good virtues of the $3^{\rm d}$ REPRESENTATION above, and where we have now a very good control of the accumulation pattern of the double point set $M_2 (f) \subset X^2$; see here Figures 2 and 3 for an illustration. Let us say that the accumulation pattern in question has sufficiently mild pathologies, so that we can proceed with our constructions. And, generically, the $M_2(f) \subset X^2$ is not a closed subset.

\smallskip

With this, starting from $X^2 \xrightarrow{ \ f \ } \widetilde M^3 (G)$, we want to build up a certain high-dimensional cell-complex $S_u (\widetilde M^3 (G))$ and one of our key steps will be to show that $S_u (\widetilde M^3 (G))$ is GSC (Theorem 8). And, in this really very fast and schematical road-map, we will not make the distinctions $S_u / S'_u$, nor $S_b / S'_b$.

\medskip

4) As a first step towards $S_u (\widetilde M^3 (G))$ we need something like a $3^{\rm d}$ thickened version but in their naive versions, neither $X^2$ and even less $fX^2$ are locally-finite. More sophisticated versions are needed (and these will certainly not be regular neighbourhoods of our objects). So, one has to drill Holes, for the present purpose of local finiteness these are closed subsets, and for later further purposes one also drills Holes which are open cells too. This way we get a locally finite $\Theta^3 (fX^2 - H ({\rm holes}))$, but the Holes will have to be compensated by 2-handles $D^2(H)$, without loosing our local finiteness. Finally we get a redefined locally finite $\Theta^3 (fX^2)$
$$
\Theta^3 (fX^2) (\mbox{now properly redefined}) = \Theta^3 (fX^2-H) + \sum_H D^2 (H) \, .
$$
So far everything is equivariant.

\medskip

5) But $\Theta^3 (fX^2-H)$ is still a singular space, with mortal singularities, the immortal ones are gone. Our singularities are all of the undrawable kind. So, in order to get smooth we go to a $4^{\rm d}$ regular neighbourhood
$$
\Theta^4 (\Theta^3 (fX^2-H) , {\mathcal R}) \, ,
$$
where ${\mathcal R}$ is a desingularization (which is not unique). Because of ${\mathcal R}$, equivariance gets lost now, while we gain smothness. By multiplying with a factor $B^N$, $N$ high, we wash out the ${\mathcal R}$-dependence and get back equivariance. This gets us to
$$
S_u (\widetilde M^3 (G)-H) \equiv \Theta^4 (\Theta^3 (fX^2-H), {\mathcal R}) \times B^N
$$
an $(N+4)$-dimensional smooth manifold. For the holes ``$-H$'' there is a system of compensating handles of index $\lambda = 2$, call this $\underset{H}{\sum} D^2 (H)$. So, finally we get our desired
$$
S_u (\widetilde M^3 (G)) = S_u (\widetilde M^3 (G)-H) + \sum_H D^2 (H) \, ,
$$
where the handle attachments are simpler (not necessarily on the boundary) and we do not correct this by thickening further. Our $S_u (\widetilde M^3 (G))$ is only a cell-complex and not a smooth manifold.

\medskip

6) Very important fact: Our $S_u$ is actually a functor, with good properties of localization and glueings. So, it makes sense downstairs too
$$
S_u ( M^3 (G)) = S_u ( M^3 (G)-H) + \sum_H D^2(H) \, ,
$$
with $S_u ( M^3 (G)-H)$ a smooth $(N+4)$-dimensional manifold. But the infinitely foamy object $S_u ( M^3 (G)-H)$ is certainly not compact. (And, would it be, our Theorem 3 would be proved already by now. But much more work is necessary.) Since everything is equivariant ant functorial, we get
$$
S_u (\widetilde M^3 (G)) = S_u ( M^3 (G))^{\sim}
$$
and, very importantly, the directly defined $S_u ( M^3 (G))$ is the quotient of $S_u (\widetilde M^3 (G))$ by $G$ (we have here two equivalent definitions).

\medskip

7) So far, our construction follows the zipping process $X^2 \Rightarrow fX^2$, which is an infinite quotient space projection. In order to get our desired GSC property, we want to replace the infinite sequence of {\bf quotient maps} of the zipping by an infinite sequence of {\bf inclusion maps}, each of which should be GSC-preserving. Here is how this goes, roughly. We start with
$$
\Theta^4 ((\Theta^3 (X^2) - H) , {\mathcal R}) \times B^N - {\rm DITCHES} \leqno (*)
$$
with an infinite system of ditches drilled in the {\bf additional dimensions} $B^N$; see here Figure 8 for an illustration. The geometric realization of the zipping consists then, essentially, in an infinite sequence of additions starting with $(*)$, a process of partial ditch filling. Each elementary step here is GSC-preserving and the whole infinite sequence of inclusion maps, mimicks the infinite sequence of quotient maps of the zipping. With this very controlled process of partial ditch-filling (and we have to proceed the way we do so as not to fall into a trap of non-metrizable mess) we create a smooth $(N+4)$-manifold $S_b (\widetilde M^3 (G) - H)$ which is such that
$$
S_b (\widetilde M^3 (G)) \equiv S_b (\widetilde M^3 (G)-H) + \sum_H D^2 (H)
$$
is GSC. The $S_b$ is again a functor and we are equivariant.

\medskip

8) In order to be able to connect $S_u$ and $S_b$ we move now downstairs at the level of $M^3 (G)$. We have an easy diffeomorphism $\eta$ entering in the following homotopy-commutative diagram
$$
\xymatrix{
S_b ( M^3 (G)-H) \ar[rr]^{\eta} &&S_u ( M^3 (G)-H) \\
&\underset{H}{\sum} \partial D^2(H) \ar[ru]_{\alpha} \ar[lu]^{\beta}
} \leqno (**)
$$
where $\beta , \alpha$ are the respective 2-handle attachments. For each $C(H) = \partial D^2(H)$ the length of the closed loop
$$
\Lambda (H) \equiv \alpha \, C(H) \underset{\gamma}{\bullet} (\eta \, \beta \, C(H))^{-1}
$$
where $\gamma$ is a connecting arc chosen such that $\Lambda (H)$ be null-homotopic, and defined using the zipping flow. The uniform boundedness of the zipping length (Theorem 4 in Part I) makes that the length of $\Lambda (H)$ is uniformly bounded. Using this fact an also the {\bf compactness} of $M^3 (G)$ we can show that $(**)$ commutes up to PROPER homotopy. It follows that, downstairs, we have
$$
S_u ( M^3 (G)) \underset{\rm DIFF}{=} S_b ( M^3 (G))
$$
and then, by functoriality we have the same at level $\widetilde M^3 (G)$. This proves Theorem 8.

\medskip

9) Finally, we want to get to the implication
$$
S_u (\widetilde M^3 (G)) \ \mbox{is GSC} \Longrightarrow \mbox{$G$ is QSF},
$$
i.e. to our Theorem 9.

\smallskip

Our $\Theta^3 (fX^2)$ is an infinitely foamy object and one changes it to a new object $\Theta^3 (fX^2)({\rm new})$ with that foaminess deleted (Lemma 16). We have a new version
$$
S_u (\widetilde M^3) = \Theta^4 (\Theta^3 (fX^2)({\rm new}),{\mathcal R}) \times B^N
$$
which continues to be equivariant and GSC. From here, proceeding like in the stabilization Lemma 5 from Part I of the survey, we get that $\Theta^4 (\Theta^3 (fX^2)({\rm new}),{\mathcal R})$ is Dehn-exhaustible and, next, that $\Theta^3 (fX^2)({\rm new})$ is itself Dehn-exhaustible.

\smallskip

This $\Theta^3 (fX^2)({\rm new})$ comes with a free $G$-action, but it is NOT co-compact. But it collapses equivariantly on a subspace $\Theta^3$ (co-compact) with a free co-compact $G$-action.

\smallskip

The proof of Theorem 9 is clinched once we prove the implication (Lemma 20)
$$
\Theta^3 (fX^2)({\rm new}) \ \mbox{Dehn-exhaustible} \Longrightarrow \Theta^3 (\mbox{co-compact}) \ \mbox{is QSF}.
$$


\newpage

\begin{thebibliography}{99}
\bibitem{1} S. Brick and M. Mihalik, The QSF property for groups and spaces, {\it Math. Z.} {\bf 200}, pp.~207-217 (1995).
\bibitem{2} J. Stallings, Brick's quasi-simple filtration for groups and 3-manifolds, {\it Geom. Group Theory} {\bf 1}, London Math. Soc. pp.~188-203 (1993).
\bibitem{3} S. Gersten and J. Stallings, Casson's idea about 3-manifolds whose universal cover is $R^3$, {\it Intern. J. of Algebra Comput.} {\bf 1}, pp.~395-406 (1991).
\bibitem{4} J. Morgan et G. Tian, Ricci Flow and the Poincar\'e Conjecture, {\it AMS, Clay Mathematical Institute} (2007).
\bibitem{5} G. Perelman, The entropy formla for the Ricci flow and its geometric applications, ArXiv:math.D6/0303109 (2002).
\bibitem{6} G. Perelman, Ricci flow with surgery on three-manifolds, ArXiv:math.D6/0303109 (2003).
\bibitem{7} G. Perelman, Finite extinction time for the solutions of the Ricci flow on certain three-manifolds, ArXiv:math.D6/0307245 (2003).
\bibitem{8} L. Bessi\`eres, G. Besson, M. Boileau, S. Maillot and J. Porti, Geometrization of $3$-manifold, {\it European Math. Soc. Tracts in Math.} {\bf 13} (2010).
\bibitem{9} V. Po\'enaru, On the equivalence relation forced by the singularities of a non degenerate simplicial map, {\it Duke Math. J.} {\bf 63}, n$^{\rm o}$~2, pp.~421-429 (1991).
\bibitem{10} V. Po\'enaru, Geometric simple connectivity and low-dimensional Topology, {\it Proc. Steklov Inst. of Math.} {\bf 247}, pp.~195-208 (2004).
\bibitem{11} V. Po\'enaru, A glimpse into the problems of the fourth dimension, preprint (2016). To appear in the Proceedings of the Conference {\it Geometry in History}, Strasbourg (2015).
\bibitem{12} V. Po\'enaru and C. Tanasi, Some remarks on geometric simple connectivity, {\it Acta Math. Hungarica} {\bf 81}, pp.~1-12 (1998).
\bibitem{13} D. Otera, V. Po\'enaru and C. Tanasi, On geometric simple connectivity, {\it Bull. Math. Soc. Sci. Math. Roumanie} {\bf 53(101)}, n$^{\rm o}$~2, pp.~157-176 (2010).
\bibitem{14} L. Funar and D.E. Otera, On the WGSC and QSF tameness conditions for finitely presented groups, ArXiv:math.GT/0610936v1 (October 2006), {\it Groups, Geometry and Dynamics} {\bf 4} n$^{\rm o}$~3, pp.~549-596 (2010).
\bibitem{15} M.W. Davis, Groups generated by reflections and aspherical manifolds not covered by Euclidean spaces, {\it Amer. Math.} {\bf 117}, pp.~293-324 (1983).
\bibitem{16} L. Bessi\`eres, La conjecture de Poincar\'e: la preuve de R. Hamilton et G. Perelman, {\it Gazette des Math.}, n$^{\rm o}$~106, pp.~7-35 (2005).
\bibitem{17} G. Besson, Une nouvelle approche de la topologie de dimension 3, d'apr\`es R. Hamilton et G. Perelman, {\it S\'eminaire Bourbaki}, $57^{\mbox{\footnotesize \`eme}}$ ann\'ee, n$^{\rm o}$~947 (2005).
\bibitem{18} V. Po\'enaru, Almost convex groups, Lipschitz combing, and $\pi_1^{\infty}$ for universal covering spaces of 3-manifolds, {\it J. Diff. Geom.} {\bf 35}, pp.~103-130 (1992).
\bibitem{19} V. Po\'enaru, Geometry \`a la Gromov for the fundamental group of a closed 3-manifold $M^3$ and the simple connectivity at infinity of $\widetilde M^3$, {\it Topology} {\bf 33}, n$^{\rm o}$~1, pp.~181-196 (1994).
\bibitem{20} V. Po\'enaru, $\pi_1^{\infty}$ and simple homotopy type in dimension 3, {\it Contemporary Math. AMS} {\bf 238}, pp.~1-28 (1999).
\bibitem{21} V. Po\'enaru, Equivariant, locally finite inverse representations with uniformly bounded zipping length for arbitrary finitely presented groups, {\it Geom. Dedicata} {\bf 167}, pp.~91-121 (2013). This is [I] of the Trilogy.
\bibitem{22} V. Po\'enaru, Geometric simple connectivity and finitely presented groups, March 2014, ArXiv/404.4283 [Math.GT]. This is [II] of the Trilogy.
\bibitem{23} V. Po\'enaru, All finitely presented groups are QSF, April 2015, ArXiv/409.7325 [Math.GT]. This is [III] of the Trilogy.
\bibitem{24} C.D. Papakyriakopoulos, On Dehn's lemma and the asphericity of knots, {\it Ann. of Math.} {\bf 66}, pp.~1-26 (1957).
\bibitem{25} J. Stallings, Group theory and Three-Dimensional Manifolds, {\it Yale Math. Monographs} {\bf 4} (1971).
\bibitem{26} V. Po\'enaru,  Killing handles of index one stably and $\pi_1^{\infty}$, {\it Duke J. Math.} {\bf 63}, n$^{\rm o}$~2, pp.~431-447 (1991).
\bibitem{27} V. Po\'enaru, The collapsible pseudo-spine representation theorem, {\it Topology} {\bf 31}, n$^{\rm o}$~3, pp.~625-636 (1992).
\bibitem{28} D. Gabai, Valentin Po\'enaru's Program for the Poincar\'e Conjecture, in the volume {\it Geometry Topology and Physics for Raoul Bott} (ed. by S.T. Yau), {\it International Press}, pp.~139-169 (1994).
\bibitem{29} V. Po\'enaru and C. Tanasi, Equivariant, almost-arborescent representations of open simply-connected 3-manifolds; a finiteness result, {\it Memoirs of the AMS} {\bf 169}, n$^{\rm o}$~800 (2004).
\bibitem{30} V. Po\'enaru and C. Tanasi,  Representations of the Whitehead manifold ${\rm Wh}^3$ and Julia sets, {\it Ann. Toulouse}, vol. IV, n$^{\rm o}$~3, pp.~665-694 (1995).
\bibitem{31} J.W. Cannon, Almost convex groups, {\it Geom. Dedicata} {\bf 22}, pp.~197-210 (1987).
\bibitem{32} M. Gromov, Infinite groups as geometric objects, {\it Proc. ICM, Warszewa}, pp.~385-392 (1983).
\bibitem{33} M. Gromov, Hyperbolic groups in Essays in group theory (S.M. Gersten, ed.), {\it Math. Sci. Res. Inst. Publ.} {\bf 8}, Springer, Berlin, pp.~75-263 (1987).
\bibitem{34} D. Sullivan and W. Thurston, Manifolds with canonical coordinate charts, some examples, {\it Enseignement Math.} {\bf 29}, pp.~15-25 (1983).
\bibitem{35} M. Freedman, The topology of four dimensional manifolds, {\it Jour. of Diff. Geometry} {\bf 68}, pp.~357-457 (1982).
\bibitem{36} M. Freedman and F. Quinn, Topology of 4-manifolds, Princeton University Press (1990).
\bibitem{37} L. Guillou and A. Marin, A la recherche de la topologie perdue, {\it Progress in Math.} {\bf 62}, Birkh\"auser (1986).
\bibitem{38} W. Thurston, The geometry and Topology of 3-manifolds, Preprint Princeton University (1979).
\bibitem{39} J. Morgan and G. Tian, The Geometrization Conjecture, {\it AMS, Clay Mathematical Institute} (2014).
\bibitem{40} D. Otera and V. Po\'enaru, ``Easy'' Representations and the QSF property for groups, {\it Bull. Belg. Math. Soc. Simon Stevin} {\bf 19}, n$^{\rm o}$~3, pp.~385-398 (2012).
\bibitem{41} V. Po\'enaru and C. Tanasi, On the Handles of index one of the product of an open simply-connected 3-manifold with a high-dimensional ball, {\it Suppl. degli Rendiconti del Circolo Matematico Palermo Serie II}, n$^{\rm o}$~62, pp.~1-135 (2000).
\bibitem{42} J. Stallings, The piece-wise linear structure of Euclidean space, {\it Proc. Cambridge Phil. Soc.} {\bf 58}, pp.~481-488 (1962).
\bibitem{43} A. Shapiro and J.H.C. Whitehead, A proof and extension of Dehn's Lemma, {\it Bull. AMS} {\bf 64}, pp.~174-178 (1958).
\bibitem{44} J. Stallings, On torsion free groups with infinitely many ends, {\it Ann. of Math.} {\bf 88}, pp.~312-334 (1968).
\bibitem{45} M. Gromov, Random walk in random groups, {\it Geom. Funct. Anal.} {\bf 3}, n$^{\rm o}$~1, pp.~73-146 (2003).
\bibitem{46} E. Ghys, Groupes al\'eatoires [d'apr\`es Misha Gromov], {\it S\'eminaire Bourbaki}, Exp.~{\bf 916} (2003).
\bibitem{47} E. Ghys, P. de la Harpe et. al., Sur les groupes hyperboliques d'apr\`es Mikhael Gromov, {\it Birkh\"auser, Progress in Math.} (1990).
\bibitem{48} M. Coornaert, T. Delzant and A. Papadopoulos, G\'eom\'etrie et th\'eorie des groupes, les groupes hyperboliques de Gromov, {\it S.L.N.} {\bf 1441} (1990).
\bibitem{49} M. Kapovich, A note in the Poenaru condition, {\it J. Group Theory} {\bf 5}, n$^{\rm o}$~1, pp.~119-127 (2002).
\bibitem{50} V. Po\'enaru, Discrete symmetry with compact fundamental domain, and geometric simple connectivity, {\it Preprint Univ. Paris-Sud Orsay} 2007-16 (2007), http://ArXiv.org/als/0711.3579.
\bibitem{51} S. Maillot, Ricci flow, scalar curvature and the Poincar\'e Conjecture, in {\it Geometry, topology, quantum field theory and cosmology}, Hermann Paris (2009) pp.~97-116.
\bibitem{52} M. Aschenbrenner, Stefan Friedle, Henry Wilton, 3-Manifold Groups, {\it EMS Series of lectures in Mathematics}, 2015.
\bibitem{53} I. Agol, The virtual Haken Conjecture, {\it Doc. Math.} {\bf 18} (2013), pp.~1045-1087.
\bibitem{54} D. Wise, The structure of groups with a quasi-convex hierarchy, pp.~1-200, Preprint 2009.
\end{thebibliography}
\end{document}